\documentclass[11pt,final]{article} 

\usepackage{amsfonts,amsthm,amsmath,amssymb, mathtools}
\usepackage{dsfont}
\usepackage{graphicx,color}
\usepackage{hyperref}
\usepackage{a4wide}
\usepackage{multirow}
\usepackage{cleveref}
\usepackage{enumitem}
\usepackage{tikz-cd}

\usepackage{mathabx}

\newtheorem{proposition}{Proposition}[section]
\newtheorem{theorem}[proposition]{Theorem}
\newtheorem{lemma}[proposition]{Lemma}
\newtheorem{corollary}[proposition]{Corollary}
\newtheorem{definition}[proposition]{Definition}

\newtheorem{problem}[proposition]{Problem}

\newenvironment{proofof}[1]{\smallskip\noindent{\textbf{Proof~of~#1.}}%
  \hspace{1pt}}{\hspace{-5pt}{\nobreak\quad\nobreak\hfill\nobreak%
    $\square$\vspace{2pt}\par}\smallskip\goodbreak}

\numberwithin{equation}{section}
\numberwithin{figure}{section}
\numberwithin{table}{section}

\allowdisplaybreaks[4]

\renewcommand{\phi}{\varphi}
\renewcommand{\theta}{\vartheta}
\renewcommand{\epsilon}{\varepsilon}
\renewcommand{\L}[1]{\mathbf{L^#1}}
\newcommand{\LL}[1]{\mathbf{L^#1}}
\renewcommand{\d}[1]{\mathinner{\mathrm{d}{#1}}}
\newcommand{\Lloc}[1]{\mathbf{L^{#1}_{loc}}}

\newcommand{\C}[1]{\mathbf{C^{#1}}}

\newcommand{\Cc}[1]{\mathbf{C_c^{#1}}}

\newcommand{\W}[2]{\mathbf{W^{#1,#2}}}
\newcommand{\BV}{\mathbf{BV}}

\newcommand{\modulo}[1]{{\left|#1\right|}}
\newcommand{\norma}[1]{{\left\|#1\right\|}}
\newcommand{\reali}{{\mathbb{R}}}
\newcommand{\naturali}{{\mathbb{N}}}

\newcommand{\tv}{\mathop\mathrm{TV}}

\newcommand{\esssup}{\mathop\mathrm{ess~sup}}

\newcommand{\spt}{\mathop\mathrm{spt}}
\newcommand{\Id}{\mathop{\mathbf{Id}}}

\newcommand{\caratt}[1]{\chi_{\strut#1}}

\newcommand{\eps}{\varepsilon}

\begin{document}

\title{Non Local Hyperbolic Dynamics of Clusters}

\author{R.M.~Colombo$^1$ \and M.~Garavello$^2$}

\maketitle

\footnotetext[1]{Unit\`a INdAM \& Dipartimento di Ingegneria
  dell'Informazione, Universit\`a di Brescia, Italy.\hfill\\
  \texttt{rinaldo.colombo@unibs.it}}

\footnotetext[2]{Dipartimento di Matematica e Applicazioni,
  Università di Milano-Bicocca, Italy.\hfill\\
  \texttt{mauro.garavello@unimib.it}}

\begin{abstract}

  \noindent The formation, movement and gluing of clusters can be
  described through a system of non local balance laws. Here, the well
  posedness of this system is obtained, as well as various stability
  estimates. Remarkably, qualitative properties of the solutions are
  proved, providing information on stationary solutions and on the
  propagation speed. In some cases, fragmentation leads to clusters
  developing independently.

  Moreover, these equations may serve as an encryption/decryption
  tool. This poses new analytical problems and asks for improved
  numerical methods.

  \medskip


  \medskip

  \noindent\textbf{MSC~2020:} Conservation Laws; Macroscopic Modeling
  of Clusters Dynamics; Pattern Formation; Hyperbolic
  Encryption/Decryption Flows
\end{abstract}


\section{Introduction}
\label{sec:introduction}

In various real-world phenomena, cluster dynamics results from the
interaction among the \emph{``individuals''} composing the system. In
many cases, fragmentation may lead to different, separate
\emph{clusters}, or \emph{patterns}, having different sizes according
to the amount of individuals focusing into them. An obvious example is
opinion dynamics, where the polarization into a single opinion,
i.e.~consensus, should not be considered as the standard type of
evolution. (Here and in what follows we adopt the terminology
in~\cite[\S~II]{Castellano2009591}).

From a macroscopic point of view, we describe the \emph{population} at
hand through its density $\rho$, which is a function of both time $t$,
varying in $\reali$, and of some sort of \emph{position} $x$, varying
in a suitable $\reali^n$. As long as individuals neither appear nor
disappear, $\rho$ solves a \emph{continuity equation}, i.e., a
\emph{conservation law}~\cite[\S~1.4]{MR3468916} of the form
\begin{equation}
  \label{eq:51}
  \partial_t \rho + \nabla \cdot (\rho\,v) = 0
\end{equation}
for a suitable velocity field $v$ in $\reali^n$.

In fluid dynamics, basic physical principles justify the use of
further conservation (or balance) laws, such as that of linear
momentum or that of energy. Social dynamics lacks these basic
principles. At the same time, striving for simplicity induces to
complete the model with an ansatz specifying how $v$ depends on
$\rho$, thereby closing~\eqref{eq:51}.  Hence, in general, the agent
(or the opinion/individual) $p$, that at time $t_o$ is at $x_o$, moves
(or evolves) according to
\begin{equation}
  \label{eq:53}
  \left\{
    \begin{array}{l}
      \dot p = v\left( t, p (t), \rho (t)\right)
      \\
      p (t_o) = x_o
    \end{array}
  \right.
\end{equation}
where the dependence of $v$ on $\rho (t)$ is of a \emph{functional}
nature. Indeed, in many cases, we expect $p$ to react to the values of
$\rho (t,x)$ for all $x$ in a neighborhood of $p$ and not only to the
value $\rho\left(t, p (t)\right)$ computed exclusively at
$\left(t, p (t)\right)$.

In the dynamics of $p (t)$ the relevance of the value of $\rho (t,x)$
at some $x \in \reali^n$ realistically depends on the distance
$\norma{x-p (t)}$ or, more generally, on the difference $x-p
(t)$. This naturally leads to postulate that
$v\left(t, p (t), \rho\right)$ depends on $\rho$ through some sort of
weighted average
$\int_{\reali^n} \eta \left(x-p (t)\right) \, \rho (t,x) \d{x}$, where
$\eta\left(x-p (t)\right)$ is the (non necessarily positive) weight of
$\rho(t, x)$ in the evolution of $p (t)$.

Whenever $\rho (t,x)$ is constant in $x$ on all $\reali^n$, it is
reasonable to expect that $v$ in~\eqref{eq:53} vanishes. In other
words, the \emph{variations} of
$\int_{\reali^n} \eta \left(x-p (t)\right) \, \rho (t,x) \d{x}$ appear
to be crucial for the dynamics of $p (t)$. Thus, we consider the model
\begin{equation}
  \label{eq:52}
  \partial_t \rho
  + \nabla \cdot \left(\rho \; V \left(\nabla (\rho * \eta)\right)\right)
  = 0
\end{equation}
where $V\colon \reali^n \to \reali^n$ is a given function and, as
usual,
$(\rho * \eta) (x) = \int_{\reali^n} \eta (\xi-x) \, \rho (t,\xi)
\d\xi$. As long as $\rho$ and $\eta$ are sufficiently smooth, the
elementary properties of the convolution product ensure that
$\nabla (\rho*\eta) = (\nabla\rho) * \eta$, meaning that the gradient
of the average $\rho*\eta$ is the average of the gradient of $\rho$.

A straightforward extension of~\eqref{eq:52} to the case of several
interacting populations is then
\begin{equation}
  \label{eq:1}
  \partial_t \rho_i
  +
  \nabla \cdot \left( \rho_i \; V_i (\nabla \rho * \eta)\right) = 0
  \qquad\qquad
  \begin{array}{l}
    (t,x)
    \in
    \reali \times \reali^n
    \\
    \rho_i = \rho_i (t,x)
    \in
    \reali
    \\
    i \in \{1, \ldots, m\}\,.
  \end{array}
\end{equation}
Under two rather simple assumptions on $V$ and $\eta$,
see~\ref{item:7} and~\ref{item:6} below, we establish the well
posedness of~\eqref{eq:1}, provide various stability estimates and
prove qualitative properties of the solutions. In particular, we
present conditions ensuring that symmetries in the initial datum
persist in the solution, a large set of stationary solutions is
exhibited, estimates on the propagation speed of the initial data are
provided as well as conditions ensuring the fragmentation of the
solution.

A peculiar property of~\eqref{eq:1} is its reversibility in time. In
spite of its ability to describe fragmentation, polarization or
consensus, \eqref{eq:1} is well posed also going backward in time. As
a consequence, for instance, it can be used to encrypt and decrypt
scalar signals (where, say, $n=1$) or images (where $n=2$ and $m$ can
be thought as the number of colors). New problems arise from this
possibility. At the analytic level, one may ask if this encryption
method is \emph{``unbreakable''}, leading to the need of investigating
the properties of the group generated by~\eqref{eq:1}. At the
numerical level, the need of a \emph{reversible} algorithm arises.
Refer to \S~\ref{subsec:encrypt-decrypt-imag} for further discussions.

Once the macroscopic evolution $t \mapsto \rho (t)$ is known, the
dynamics of every single member $p_i$ of the $i$-th population is
given by~\eqref{eq:53}, which in the present multipopulation case
reads
\begin{equation}
  \label{eq:54}
  \left\{
    \begin{array}{l}
      \dot p_i = V_i\left(\nabla \left(\rho (t)*\eta\right)(p_i)\right)
      \\
      p_i (t_o) = x_o \,.
    \end{array}
  \right.
\end{equation}
The general well posedness of~\eqref{eq:54}, which we leave to a
future work, seems reachable extending the results in~\cite{MR2006201,
  MR4739982}. On the contrary, we note that the reverse connection,
i.e., a general rigorous derivation of~\eqref{eq:1} from a microscopic
model, is apparently still an open problem if $n\geq 2$. For the case
$n=1$, this connection was obtained in~\cite{MR3356989, MR3721873}, a
recent result devoted to the case $n=2$ is~\cite{weissen2021density}.

A milestone in the modeling of flocks is the Cucker-Smale
model~\cite{MR2324245}. Differently from the model presented therein,
\eqref{eq:1} is of a macroscopic nature and may describe, besides
consensus, also polarization and fragmentation, already in the case of
a single population. Moreover, with \emph{ad hoc} choices of the
function $\eta$, it can also describe interactions that are attractive
in some regions and repelling in other, for instance.

Non local equations similar to~\eqref{eq:1} have been considered in
the literature in a variety of applications. Their use dates back at
least to~\cite{MR2300313}, devoted to opinion
formation. Recall~\cite{MR4160250, MR2902155} devoted to crowd
dynamics and~\cite{zbMATH07697946} where a system of balance laws
describes flocking phenomena. A controlled diffusion equivalent to
that proved in Lemma~\ref{lem:easy} is shown in~\cite{MR4385602} in
the case of a gradient flow resulting in a system of partial
differential equations. A $1$-dimensional case of~\eqref{eq:1} has
been thoroughly studied in~\cite{MR4410035}, where solutions are
constructed as limits of solutions to interacting particle systems,
also obtaining well posedness. In~\cite{MR3414245}, non local terms
represent the competition between species for resources. The issue of
describing the effects of communications among individuals is
investigated in~\cite{MR3343585} through a non local hyperbolic
system.

Alternative to the use, here pursued, of density functions typically
in $\L1$, is the use of measures and of the various types of
Wasserstein distances. Paradigmatic in this direction is the recent
work~\cite{MR4688860}, where a general framework for population
dynamics is based on Radon measures in Polish spaces.

\medskip


The next section presents the analytical
results. Section~\ref{sec:numer-integr} is devoted to numerical
integrations. Proofs are exposed in
Section~\ref{sec:technical-proofs}, a more technical one being
deferred to the Appendix.

\section{Well Posedness and Qualitative Properties}
\label{sec:analytic-results}

In the non linear and non local equation~\eqref{eq:1},
$(\nabla \rho * \eta)$ is the $n\times m$ matrix
\begin{equation}
  \nabla \rho * \eta
  = \left[
    \begin{array}{@{}c@{}}
      \nabla \rho_1 * \eta
      \\
      \vdots
      \\
      \nabla \rho_m * \eta
    \end{array}
  \right]^\intercal
  \quad
  \begin{array}{c}
    (\nabla \rho * \eta)_{ji}
    =
    \partial_{x_j} (\rho_i * \eta)
    \mbox{ for }
    j \in \{1, \ldots, n\} \mbox{ and } i \in \{1, \ldots, m\}
    \\[6pt]
    \mbox{where } \quad
    \displaystyle
    (\rho_i (t) * \eta) (x)
    =
    \int_{\reali^n} \rho_i (t,\xi) \; \eta (x-\xi) \, \d\xi \,.
  \end{array}
\end{equation}
The extension to different kernels $\eta_{ij}$ acting on the different
populations requires, from the analytical point of view, essentially
only formal modifications, see~\cite{Claudia01, MR3057143}.

Throughout, we pose the following assumptions:
\begin{enumerate}[label=$(\mathbf{\eta}$)]
\item \label{item:6}
  $\eta \in (\C3 \cap \W31 \cap \W3\infty) (\reali^n; \reali)$.
\end{enumerate}
\begin{enumerate}[label=$(\mathbf{V})$]
\item \label{item:7}
  $V \in \C2 (\reali^{n\times m}; \reali^{n\times m})$, $V (0) = 0$
  and
  $\norma{DV}_{\W1\infty (\reali^{n\times m};\reali^{(n\times m)^2})}
  \leq L_V$, for a $L_V \in \reali_+$.
\end{enumerate}
\noindent We set $V = [V_1 \cdots V_m]$, with $V_i \in \reali^n$.

The next definition is intrinsic to~\eqref{eq:1} and independent of
the way in which any particular solution may be constructed.

\begin{definition}
  \label{def:solution}
  Fix a non trivial real interval $I$.  A map
  $\rho \in \C0 \left(I;\L1 (\reali^n; \reali^m)\right)$ is a
  \emph{solution} to~\eqref{eq:1} on $I$ if setting, for
  $i \in \{1, \ldots, m\}$ and $(t,x) \in I \times \reali^n$,
  \begin{equation}
    \label{eq:25}
    v_i (t,x)
    =
    \left[
      V_{1i} \left((\nabla \rho * \eta) (t,x)\right)
      \quad V_{2i} \left((\nabla \rho * \eta) (t,x)\right)
      \quad \cdots
      \quad V_{ni} \left((\nabla \rho * \eta) (t,x)\right)
    \right]^\intercal \,,
  \end{equation}
  for $i \in \{1, \ldots, m\}$ the component $\rho_i$ is a solution to
  \begin{equation}
    \label{eq:23}
    \partial_t \rho_i
    +
    \nabla \cdot \left( \rho_i \, v_i (t,x)\right) = 0
    \qquad (t,x) \in I \times \reali^n
  \end{equation}
  on the interval $I$.
\end{definition}
\noindent Refer to~\eqref{eq:7} and to the discussion therein for the
precise meaning of \emph{solution} to~\eqref{eq:23}.

Here follows the main result of the present paper.

\begin{theorem}
  \label{thm:cauchyProblem}
  Let~\ref{item:6} and~\ref{item:7} 
  hold. Then~\eqref{eq:1} generates a unique map
  \begin{displaymath}
    \mathcal{G} \colon \reali \times
    (\L1 \cap \L\infty \cap \BV) (\reali^n; \reali^m)
    \to
    (\L1 \cap \L\infty \cap \BV) (\reali^n; \reali^m)
  \end{displaymath}
  such that
  \begin{enumerate}[label={\bf(\arabic*)}]

  \item \label{item:33} $\mathcal{G}$ is a group, in the sense that
    \begin{displaymath}
      \mathcal{G}_0 = \Id
      \quad \mbox{ and } \quad
      \mathcal{G}_{t_1} \circ \mathcal{G}_{t_2} = \mathcal{G}_{t_1 + t_2}
    \end{displaymath}
    for every $t_1, t_2 \in \reali$.
  \item \label{item:32} For any
    $\rho_o \in (\L1 \cap \L\infty \cap \BV) (\reali^n; \reali^m)$,
    the map $t \mapsto \mathcal{G}_t \rho_o$ is the unique global
    solution to~\eqref{eq:1} with initial datum $\rho_o$ assigned at
    time $t=0$, in the sense of Definition~\ref{def:solution}.
  \item \label{item:constant_L1} For any
    $\rho_o \in (\L1 \cap \L\infty \cap \BV) (\reali^n; \reali^m)$ and
    for every $t \in \reali$,
    \begin{equation*}
      \norma{\mathcal{G}_t \rho_o}_{\L1 (\reali^n; \reali^m)}
      = \norma{\rho_o}_{\L1 (\reali^n; \reali^m)}.
    \end{equation*}

  \item \label{item:9} For any
    $\widehat\rho_o, \widecheck\rho_o \in (\L1 \cap \L\infty \cap \BV)
    (\reali^n; \reali^m)$ and for every $t \in \reali$,
    \begin{displaymath}
      \norma{\mathcal{G}_t \widehat \rho_o - \mathcal{G}_t\widecheck \rho_o}
      _{\L1 (\reali^n; \reali^m)}
      \leq
      \norma{\widecheck \rho_o - \widehat \rho_o}
      _{\L1 (\reali^n; \reali^{m})} \;
      e^{C (|t|)}
    \end{displaymath}
    where the function $C \in \C0 (\reali_+; \reali_+)$ is non
    decreasing, depends only on $\eta$, $L_V$ and on the $\L1$ and
    $\L\infty$ norms of $\widehat\rho_o$, $\widecheck\rho_o$; moreover
    $\limsup_{t\to 0+} C (t)/t$ is bounded.

  \item \label{item:14} If $\widehat{\mathcal{G}}$, respectively
    $\widecheck{\mathcal{G}}$, is generated by~\eqref{eq:1} with speed
    $\widehat V$, respectively $\widecheck V$, satisfying~\ref{item:7}
    and such that $L_V$ is an upper bound for the $\W1\infty$ norms of
    both $D\widehat V$ and $D\widecheck V$, then, for every
    $t \in \reali$ and for any
    $\rho_o \in (\L1 \cap \L\infty \cap \BV) (\reali^n; \reali^m)$
    \begin{displaymath}
      \norma{ \widehat {\mathcal G}_t \rho_o - \widecheck {\mathcal
          G}_t \rho_o}_{\L1 (\reali^n; \reali^m)}
      \leq
      C \, \modulo{t}
      \left(
        \norma{\widehat{V} - \widecheck{V}}_{\L\infty (B;\reali^{n\times m})}
        +
        \norma{D\widehat{V} - D\widecheck{V}}_{\L1 (B;\reali^{\left(n\times m\right)^2})}
      \right) e^{C\,\modulo{t}}
    \end{displaymath}
    where
    $B = B_{\reali^{n\times m}} (0, \norma{\rho_o}_{\L1 (\reali^n;
      \reali^m)} \, \norma{\nabla\eta}_{\L\infty (\reali^n;
      \reali^n)})$ and $C$ depends on various norms of $\eta$ and of
    the initial datum.
  \item \label{item:12} For $k = 0,1,2$, if
    $\rho_o \in \C{k} (\reali^n; \reali^m)$, then also
    $\rho (t) \in \C{k} (\reali^n; \reali^m)$ for all
    $t \in \reali_+$.
  \item \label{item:4} For any $i \in \{1, \ldots, m\}$, if
    $\rho_{o,i} \geq 0$ then for all
    $(t,x) \in \reali_+ \times \reali^n$, also $\rho_i (t,x) \geq 0$.
  \end{enumerate}
\end{theorem}

\noindent The proof relies on a fixed point argument and is deferred
to Section~\ref{sec:technical-proofs}. Note that the group property
stated at~\ref{item:33} implies that problem~\eqref{eq:1} is
\emph{reversible}.

Next, we provide conditions under which~\eqref{eq:1} is invariant with
respect to symmetries in $O (n)$, i.e.~orthogonal matrices.

\begin{proposition}
  \label{prop:polar}
  Let~\ref{item:7} and~\ref{item:6} hold. Fix
  $\rho_o \in (\L1\cap \L\infty \cap \BV) (\reali^n;
  \reali^m)$. Assume moreover that for a $R \in O (n)$:
  \begin{enumerate}[label=\bf(\arabic*)]
  \item \label{item:19}
    $V_i (w_1\, R, \ldots, w_m\, R) = R^{-1} \, V_i (w_1, \ldots,
    w_m)$ for $i\in \{1, \ldots, m\}$ and all
    $w_1, \ldots, w_m \in \reali^{1\times n}$.
  \item \label{item:18} $\eta (R\,x) = \eta (x)$ for all
    $x \in \reali^n$.
  \item \label{item:20} $\rho_o (R\,x) = \rho_o (x)$ for all
    $x \in \reali^n$.
  \end{enumerate}
  Then, the solution $t \mapsto \mathcal{G}_t \rho_o$ exhibited in
  Theorem~\ref{thm:cauchyProblem} is invariant with respect to $R$,
  i.e.,
  \begin{displaymath}
    \forall (t,x) \in \reali \times \reali^n \qquad
    (\mathcal{G}_t \rho_o)(R\,x) = (\mathcal{G}_t \rho_o) (x) \,.
  \end{displaymath}
\end{proposition}


On the basis of~\ref{item:7} it is immediate to provide the following
general bound on the growth of the support of the solutions
to~\eqref{eq:1}. This shows that solutions to~\eqref{eq:1} propagate
with finite speed.

\begin{lemma}
  \label{lem:easy}
  Let~\ref{item:6} and~\ref{item:7} hold. Fix
  $\rho_o \in (\L1\cap\L\infty\cap\BV) (\reali^n; \reali^m)$,
  $x_o \in \reali^n$ and $r>0$. If $\spt \rho_o \subseteq B (x_o,r)$,
  then for all $t \in \reali$
  \begin{equation}
    \label{eq:9}
    \spt \rho (t) \subseteq B \left(x_o, r + W\,\modulo{t}\right)
    \quad \mbox{ where } \quad
    W
    =
    L_V \; \norma{\nabla\eta}_{\L\infty (\reali^n; \reali)}
    \; \norma{\rho_o}_{\L1 (\reali^n; \reali^m)} \,.
  \end{equation}
\end{lemma}

As a consequence, a compactly supported initial datum yields a
solution which is compactly supported at any positive time.

The qualitative behavior of the solutions to~\eqref{eq:1} is governed
by the direction of the vectors $V_i$. Indeed, when $V_i$ and
$\nabla (\rho_i *\eta)$ have the same direction,
i.e.~$V_i \left(\nabla (\rho *\eta)\right)\cdot \nabla (\rho_i
*\eta)>0$, then solutions to~\eqref{eq:1} do not propagate, in the
sense of the following Proposition.

\begin{proposition}
  \label{prop:Caratt}
  Let~\ref{item:7}, \ref{item:6} hold and
  $\rho_o \in (\L1\cap\L\infty\cap\BV) (\reali^n; \reali^m)$. Assume
  moreover
  \begin{enumerate}[label=\bf(\arabic*)]
  \item \label{item:17} For all $w \in \reali^{n\times m}$,
    $V (w) = v \left(\norma{w}\right) \, w$ with $v$ such that
    $v (s) >0$ for all $s>0$.
  \item \label{item:21} For all $x \in \reali^n$,
    $\eta (x) = \tilde\eta \left(\norma{x}\right)$ with
    $\tilde\eta' \leq 0$ and $\spt\tilde\eta = [-\ell,\ell]$.
  \item \label{item:24} For all $i \in \{1, \ldots, m\}$,
    $\rho_{o,i} (x) \geq 0$ for a.e.~$x \in \reali^n$.
  \end{enumerate}
  Then,
  \begin{enumerate}[label=\bf(\roman*)]
  \item \label{item:22} If $C \subseteq \reali^n$ is closed and convex
    and $i \in \left\{1, \ldots, m\right\}$,
    \begin{displaymath}
      \spt \rho_{o,i}\subset C
      \implies
      \forall\, t \in \reali_+ \quad
      \spt (\mathcal{G}_t \rho)_i \subseteq C \,.
    \end{displaymath}
  \item \label{item:23} If $\spt \rho_o$ is bounded, then for all
    $t \in \reali_+$ and $i \in \left\{1, \ldots, m\right\}$,
    \begin{displaymath}
      \spt (\mathcal{G}_t \rho)_i
      \subseteq \overline{\mathop{\rm co}} \spt \rho_{o,i} \,.
    \end{displaymath}
  \end{enumerate}
\end{proposition}

\noindent Above, we adopt the usual definition of the support of a
function $\rho \in \L1 (\reali^n; \reali^m)$:
\begin{equation}
  \label{eq:41}
  \spt\rho
  =
  \reali^n \Big\backslash \bigcup_{A \in \mathcal{A}} A
  \quad \mbox{ where } \quad
  \mathcal{A}
  =
  \left\{A \subseteq \reali^n \colon A \mbox{ is open and
    } \rho = 0 \mbox{ a.e.~in } A \right\}\,,
\end{equation}
see~\cite[Proposition~IV.17]{Brezis}, and
$\overline{\mathop{\rm co}} \, S$ stands for the closed convex hull of
the set $S$.

Recall that $\ell$ is the radius of the support of $\eta$, i.e., the
maximal distance at which the non local interaction may act. We now
show that, under the key assumption~\ref{item:17} of
Proposition~\ref{prop:Caratt}, if the initial datum is supported over
balls at a distance at least $\ell$ from each other, then the
evolution of the solution on these balls proceeds independently on the
different spheres.

\begin{figure}[!h]
  \centering
  \begin{tikzpicture}[line cap=round,line join=round,x=1.cm,y=1.cm]

    \draw[] (3. , 1.) circle (0.6cm);%
    \draw[fill] (3. , 1.) circle (0.03cm); %
    \draw[dashed] (3. , 1.)  circle (1.3cm); %
    \draw[->, color=red] (3., 1.) -- (3., 1.6); %
    \node[inner sep=0, anchor=north west] at (3.1, 1.) {$x_1$}; %
    \node[inner sep=0, anchor=west] at (3.1, 1.3)
    {\color{red}$r_1$}; %

    \draw[->, color=blue] (3., 1.) -- (1.7, 1.); \node[inner sep=0,
    anchor=north] at (2.35, .9) {\color{blue}$r_1 + \ell$};


    \draw[] (6. , 3.) circle (0.8cm); %
    \draw[fill] (6. , 3.) circle (0.03cm); %
    \draw[dashed] (6. , 3.) circle (1.5cm); %
    \draw[->, color=red] (6., 3.) -- (6., 3.8); %
    \draw[->, color=blue] (6., 3.)  -- (4.5, 3.); %
    \node[inner sep=0, anchor=north] at (5.25, 2.9)
    {\color{blue}$r_2 + \ell$}; %

    \node[inner sep=0, anchor=north west] at (6.1, 3.) {$x_2$};
    \node[inner sep=0, anchor=west] at (6.1, 3.3) {\color{red}$r_2$};

    \draw[] (9. , 2.+0.4) circle (1.2cm); %
    \draw[fill] (9. , 2.+0.4) circle (0.03cm); %
    \draw[dashed] (9. , 2.+0.4) circle (1.9cm); %
    \draw[->, color=red] (9., 2.+0.4) -- (9., 3.2+0.4); %
    \draw[->, color=blue] (9., 2.+0.4) -- (7.1, 2.+0.4); %
    \node[inner sep=0, anchor=north] at (8.05, 1.9+0.4)
    {\color{blue}$r_3 + \ell$}; %

    \node[inner sep=0, anchor=north west] at (9.1, 2.+0.4) {$x_3$};
    \node[inner sep=0, anchor=west] at (9.1, 2.3+0.4)
    {\color{red}$r_3$};
  \end{tikzpicture}

  \caption{An example of a situation complying with~\eqref{eq:32} of
    Corollary~\ref{cor:clustering} in the case $k = 3$. The dashed
    circumferences may not overlap with the solid ones.}
  \label{fig:clusters}
\end{figure}
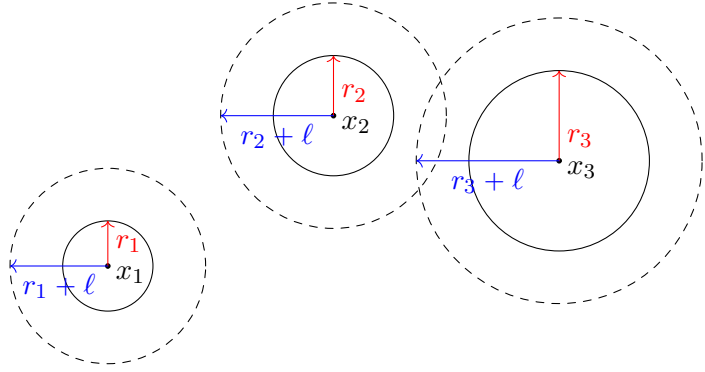

\begin{corollary}
  \label{cor:clustering}
  Under the assumptions~\ref{item:17}, \ref{item:21} and~\ref{item:24}
  in Proposition~\ref{prop:Caratt}, assume that, for a suitable
  $i \in \{1, \ldots, m\}$, there exist $x_1, \ldots, x_k\in \reali^n$
  and $r_1, \ldots, r_k \in \mathopen]0, +\infty\mathclose[$, such
  that
  \begin{equation}
    \label{eq:32}
    \spt \rho_{o,i} \subseteq \bigcup_{h=1}^k B (x_h,r_h)
    \mbox{ with }
    \norma{x_h - x_j} > r_h + r_j + \ell
    \mbox{ for all } h,j \in \{1, \ldots, k\}\,,h\neq j \,.
  \end{equation}
  Then, the solution $\mathcal{G} \rho_o$ to~\eqref{eq:1} satisfies
  for all $t \in \reali_+$
  \begin{eqnarray}
    \label{eq:28}
    \spt (\mathcal{G}_t \rho_o)_i
    & \subseteq
    & \bigcup_{h=1}^k B (x_h,r_h)
      \qquad\qquad \mbox{ and}
    \\
    \label{eq:26}
    \int_{B (x_h,r_h)} (\mathcal{G}_t \rho_o)_i (x) \d{x}
    & =
    & \int_{B (x_h,r_h)} \rho_{o,i} (x) \d{x}
      \qquad \mbox{for } h \in \{1, \ldots, k\}\,.
  \end{eqnarray}
\end{corollary}

The next proposition describes a set of stationary solutions
to~\eqref{eq:1}: they are the most natural candidates of time
asymptotic limits of solutions to~\eqref{eq:1}.

\begin{proposition}
  \label{prop:stationary}
  Let~\ref{item:7} and~\ref{item:6} hold. Fix
  $\rho_o \in (\L1\cap \L\infty\cap \BV) (\reali^n; \reali^m)$. Assume
  moreover that for fixed $r,\ell \in \reali$ with $\ell > r > 0$ and
  for $x_1, \ldots, x_k \in \reali^n$,
  \begin{enumerate}[label=\bf(\arabic*)]
  \item \label{item:30}
    $\spt \nabla \eta \subseteq B(0,\ell) \setminus B (0,r)$.
  \item \label{item:31}
    $\spt \rho_o \subseteq \bigcup_{h=1}^k B (x_h,r / 2)$ with
    $\norma{x_h - x_j} > r + \ell$ for all $h,j \in \{1, \ldots, k\}$
    with $h\neq j$.
  \end{enumerate}
  Then, the solution to~\eqref{eq:1} is stationary:
  $\mathcal{G}_t \rho_o = \rho_o$ for all $t\geq 0$.
\end{proposition}

In the case $k=1$, assumption~\ref{item:31} in
Proposition~\ref{prop:stationary} simplifies to
$\spt \rho_o \subseteq B (x_1,r / 2)$.

\section{Examples}
\label{sec:numer-integr}

The paragraphs below present numerical integrations
of~\eqref{eq:1}. They are obtained by means of the upwind
method~\cite[Formula~(4.35) in \S~4.8]{LeVequeBook2002} or of the
Lax-Friedrichs method~\cite[\S~4.6]{LeVequeBook2002}. The space mesh
is uniform while the time step is chosen adaptively to meet the
(adapted) CFL condition~\cite[\S~4.4]{LeVequeBook2002}. Whenever
$n>1$, we adopt the usual dimensional
splitting~\cite[\S~19.5]{LeVequeBook2002}. The convolutions are
approximated through integrals of functions that are piecewise
constant on the spatial mesh.

Below, initial data are often chosen to be linear combination of
characteristic functions. This choice, though particular, is balanced
by the $\L1$ stability proved in Theorem~\ref{thm:cauchyProblem}.

\subsection{Fragmentation, or Clusters Formation}
\label{subsec:cluste-formation}


Here, we exemplify the situation described in
Corollary~\ref{cor:clustering}. To this aim, choose the initial datum and the parameters as in~\eqref{eq:50}:\\
\begin{minipage}{0.7\linewidth}
  \includegraphics[width=0.5\linewidth,trim=40 0 40
  0]{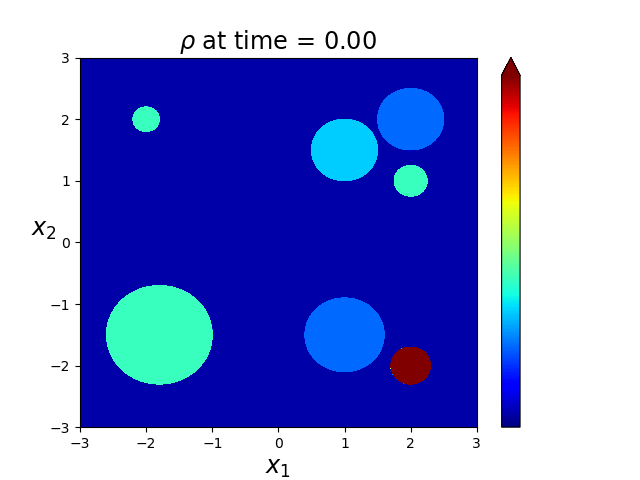}%
  \includegraphics[width=0.5\linewidth,trim=40 0 40
  0]{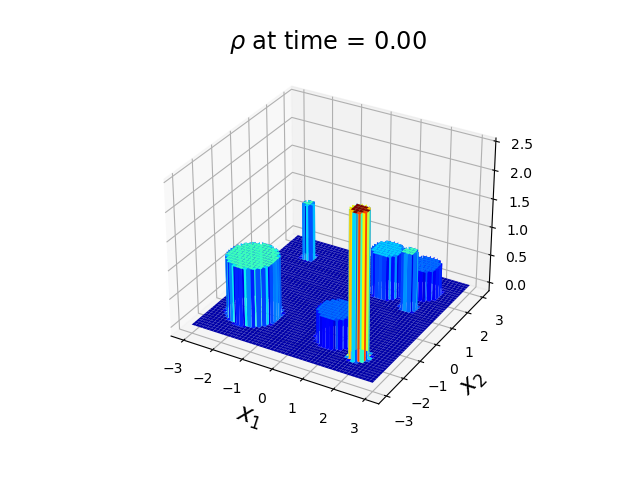}
\end{minipage}%
\begin{minipage}{0.3\linewidth}
  \begin{equation}
    \label{eq:50}
    \begin{array}{@{}r@{\,}c@{\,}l@{}}
      n
      & =
      & 2 \,,
      \\
      m
      & =
      & 1 \,,
      \\
      V (w)
      & =
      & \dfrac{w}{\sqrt{1{+}\norma{w}^2}} \,.
    \end{array}
  \end{equation}
\end{minipage}\\
More precisely, the initial datum is a linear combination of
characteristic functions of circles
$\rho_o = \sum_{\nu=1}^7 h_\nu \, \caratt{B (c_\nu,r_\nu)}$ where
\begin{equation}
  \label{eq:56}
  \begin{array}{@{}cccccccc@{}}
    \nu
    & 1
    & 2
    & 3
    & 4
    & 5
    & 6
    & 7
    \\ h_\nu
    & 1.0
    & 1.0
    & 0.5
    & 0.75
    & 1
    & 0.5
    & 2.5
    \\ c_\nu
    & (-2,2)
    & (-1.8,1.5)
    & (2,2)
    & (1, 1.5)
    & (2,1)
    & (1, -1.5)
    & (2, -2)
    \\r_\nu
    & 0.2
    & 0.8
    & 0.5
    & 0.5
    & 0.25
    & 0.6
    & 0.3
  \end{array}
\end{equation}
as shown in the two diagrams on the left in~\eqref{eq:50}. As function
$\eta$ we use\\
\begin{minipage}{0.4\linewidth}
  \includegraphics[width=\linewidth,trim=40 20 40
  5]{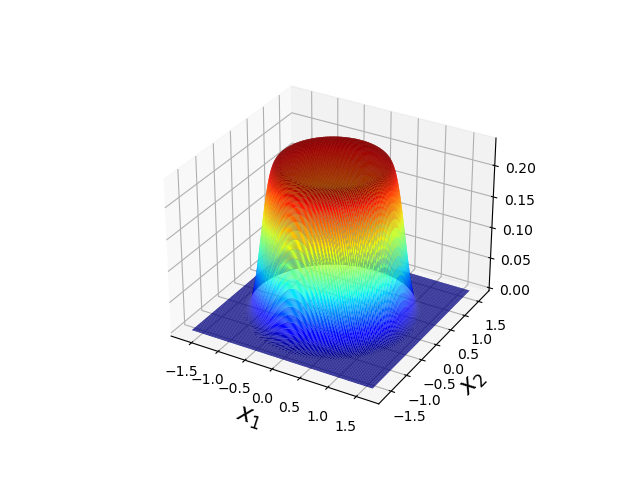}
\end{minipage}%
\begin{minipage}{0.6\linewidth}
  \begin{equation}
    \label{eq:46}
    \begin{array}{@{}rcl@{}}
      a (\xi)
      & =
      & \left\{
        \begin{array}{ll}
          0
          & \xi < r
          \\
          k \, (\xi-r)^3 \, (\ell-\xi)^3
          & \xi \in [r,\ell]
          \\
          0
          & \xi > \ell
        \end{array}
            \right.
      \\
      \eta (x)
      & =
      & \displaystyle \int_{\norma{x}}^\ell a (\xi) \, \d\xi
      \\
      k
      & \colon
      & \displaystyle \int_{\reali^2} \eta (x) \d{x} =1\,.
    \end{array}
  \end{equation}
\end{minipage}
\\ with $r = 0.5$ and $\ell = 0.8$.

\begin{figure}[!h]
  \centering \includegraphics[width=0.25\linewidth,trim=50 0 40
  0]{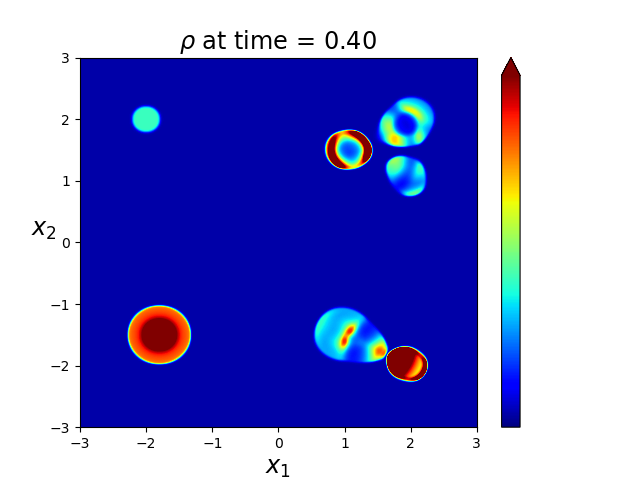}%
  \includegraphics[width=0.25\linewidth,trim=50 0 40
  0]{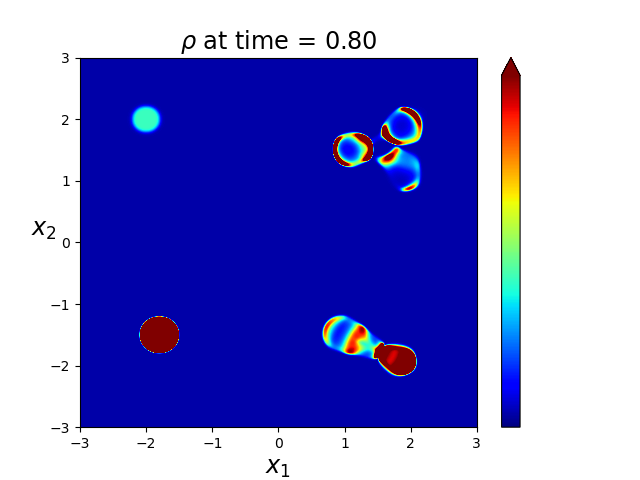}%
  \includegraphics[width=0.25\linewidth,trim=50 0 40
  0]{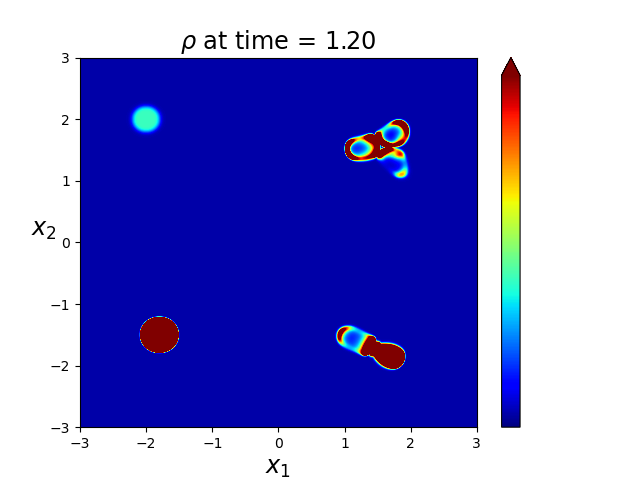}%
  \includegraphics[width=0.25\linewidth,trim=50 0 40 0]{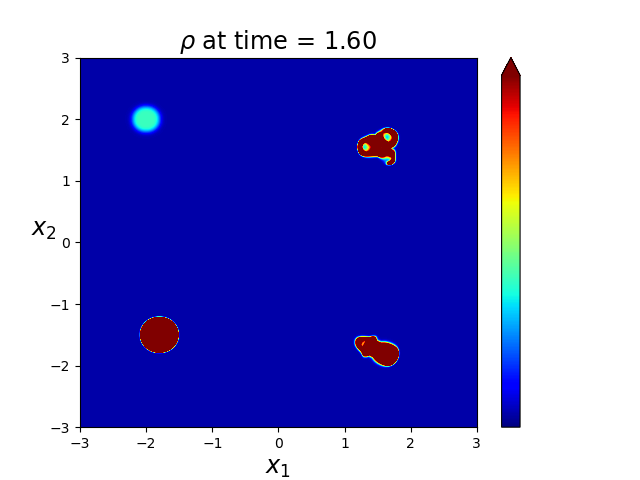}\\
  \caption{Solution to~\eqref{eq:1}--\eqref{eq:50} with $\eta$ as
    in~\eqref{eq:46}. This evolution is consistent with
    Corollary~\ref{cor:clustering}, since the $4$ parts of the
    solutions in the $4$ corners evolve independently, as well as with
    Proposition~\ref{prop:stationary}, since the top left part is
    stationary.}
  \label{fig:clPPLF}
\end{figure}

The behavior of the solution is consistent with
Corollary~\ref{cor:clustering}. Indeed, as Figure~\ref{fig:clPPLF}
shows, the $4$ clusters of characteristic functions in the $4$ corners
of the numerical domain $[-3,\, 3]^2$ evolve independently. In the top
left corner we have a stationary cluster, consistent with
Proposition~\ref{prop:stationary}. In the bottom left corner, the
initial cluster is supported in a circle too wide to be stationary and
evolves focusing in a smaller circle, its volume being conserved.

The $2$ clusters on the right consist of different initial
characteristic functions. Each of the $2$ clusters focuses in a sort
of barycenter of the initial mass.

Here, we adopted the Lax-Friedrichs scheme on the numerical domain
$[-3,3]^2$ with $7000\times7000$ points.

\subsection{On the Role of \texorpdfstring{$V$}{V}}


This paragraph is devoted to exemplify the role of $V$ in~\eqref{eq:1}. Consider the initial datum and the parameters in~\eqref{eq:55}\\
\begin{minipage}{0.7\linewidth}
  \includegraphics[width=0.5\linewidth,trim=40 0 40 0]{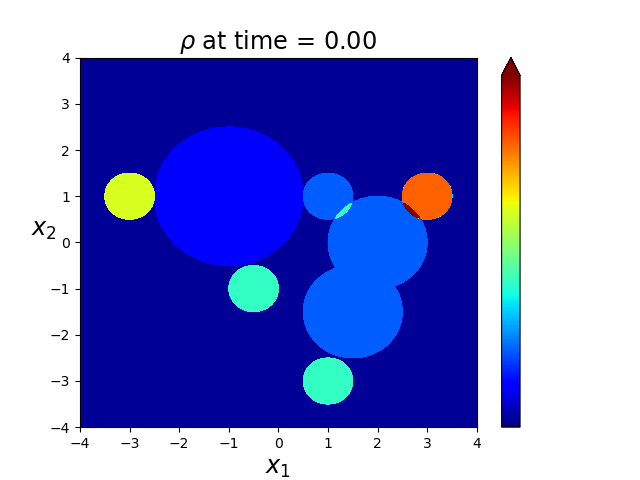}%
  \includegraphics[width=0.5\linewidth,trim=40 0 40 0]{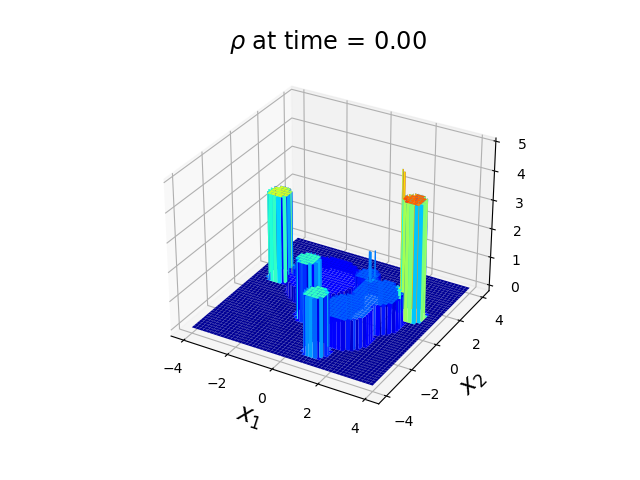}
\end{minipage}%
\begin{minipage}{0.3\linewidth}
  \begin{equation}
    \label{eq:55}
    \begin{array}{@{}r@{\,}c@{\,}l@{}}
      n
      & =
      & 2 \,,
      \\
      m
      & =
      & 1 \,,
      \\
      r
      & =
      & 0.6 \,,
      \\
      \ell
      & =
      & 1.5
    \end{array}
  \end{equation}
\end{minipage}\\
with $\eta$ as in~\eqref{eq:46} and with the following $3$ choices for
$V$:
\begin{equation}
  \label{eq:57}
  V^1 (w)
  =
  \dfrac{w}{\sqrt{1{+}\norma{w}^2}} \,,
  \qquad
  V^2 (w)
  =
  \left[
    \begin{array}{@{}cc@{}}
      0
      & 1
      \\
      -1
      & 0
    \end{array}
  \right]
  \dfrac{w}{\sqrt{1{+}\norma{w}^2}} \,,
  \qquad
  V^3 (w)
  =
  \dfrac{-w}{\sqrt{1{+}\norma{w}^2}} \,.
\end{equation}
Only the former choice $V^1$ complies with~\ref{item:17} in
Corollary~\ref{cor:clustering} and, indeed, in this case the support
of the clusters concentrate. In the other $2$ cases, some mass exits
the numerical domain.
\begin{figure}[!h]
  \includegraphics[width=0.33\linewidth,trim=40 0 40
  0]{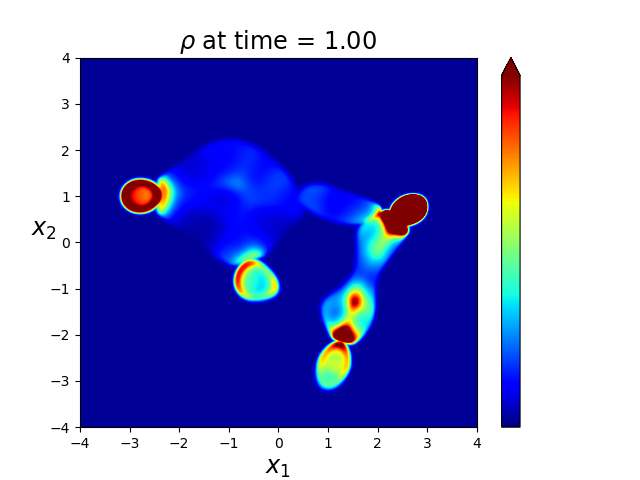}%
  \includegraphics[width=0.33\linewidth,trim=40 0 40
  0]{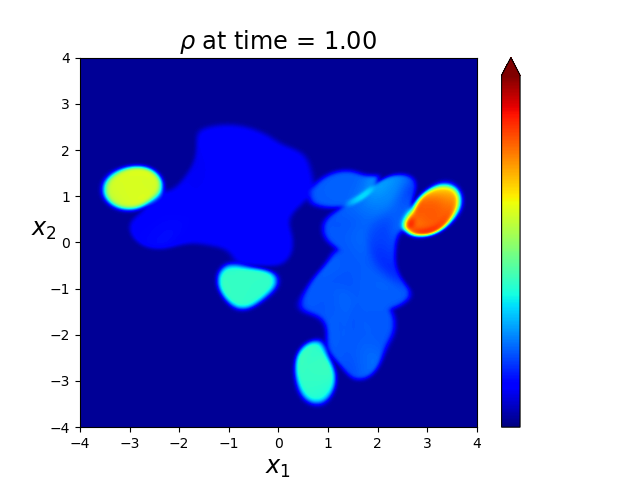}%
  \includegraphics[width=0.33\linewidth,trim=40 0 40 0]{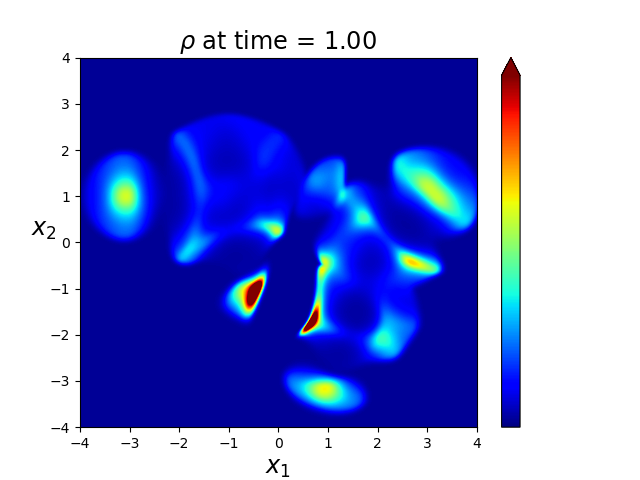}\\
  \includegraphics[width=0.33\linewidth,trim=40 0 40
  0]{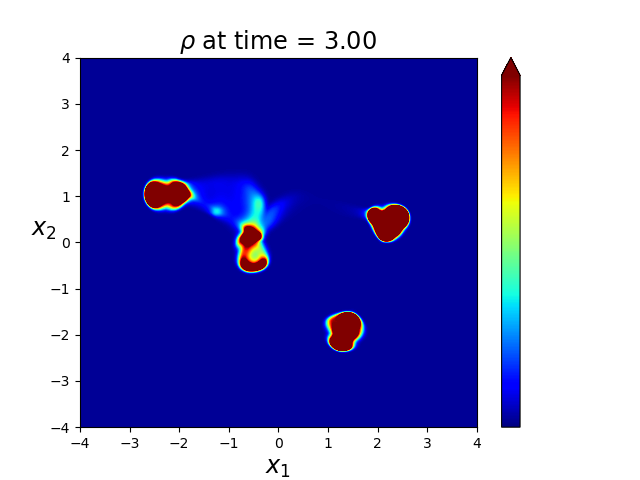}%
  \includegraphics[width=0.33\linewidth,trim=40 0 40
  0]{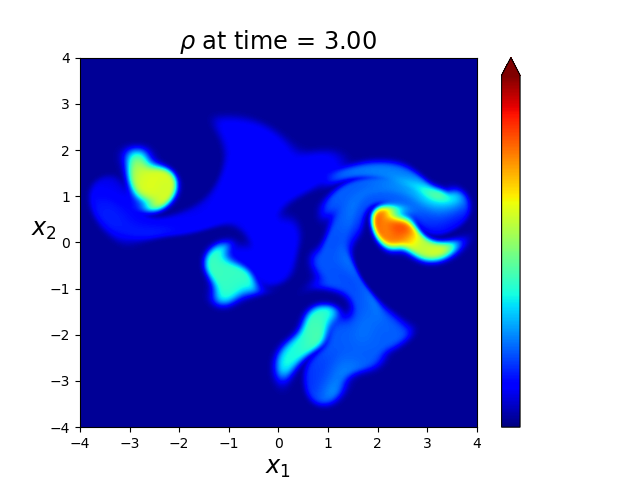}%
  \includegraphics[width=0.33\linewidth,trim=40 0 40 0]{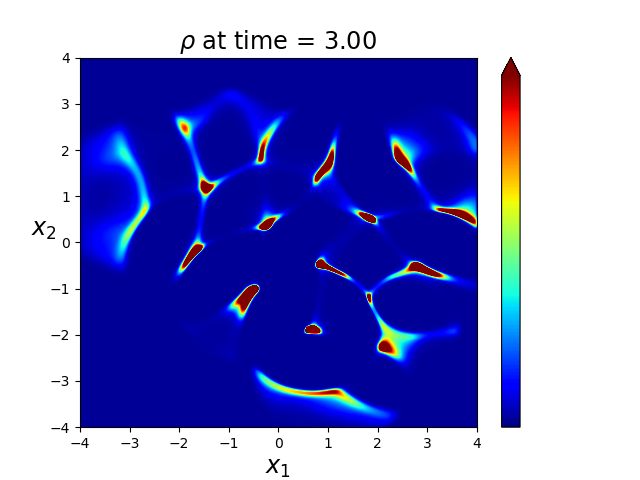}\\
  \caption{Evolution of the solutions
    to~\eqref{eq:1}--\eqref{eq:55}--\eqref{eq:46}. On the left,
    respectively middle and right column, the vector field is $V^1$,
    respectively $V^2$ and $V^3$ in~\eqref{eq:57}. Only in the
    solution in the left column is the total mass conserved. In the
    middle and right cases, some mass exits the numerical domain.}
  \label{fig:bella2}
\end{figure}
Note that $\eta$ as defined in~\eqref{eq:46} satisfies~\ref{item:6}
and~\ref{item:30} in Proposition~\ref{prop:stationary}. All $V^1$,
$V^2$ and $V^3$ comply with~\ref{item:7}. Moreover, the pattern
formation resulting in all $3$ cases apparently leads asymptotically
to separate clusters, see Figure~\ref{fig:bella2}, corresponding to
the stationary solutions exhibited in
Proposition~\ref{prop:stationary}.

All $3$ integrations were obtained by means of Lax-Friedrichs method
on the numerical domain $[-4,4]^2$ with $3000\times3000$ mesh points.

\subsection{A Stationary Solution}
\label{subsec:stationary-solution}


As an exemplification of the stationary solutions exhibited in
Proposition~\ref{prop:stationary}, consider the following setting:
\begin{equation}
  \label{eq:45}
  \begin{array}{r@{\,}c@{\,}l}
    n
    & =
    & 2 \,,
    \\
    m
    & =
    & 1 \,,
  \end{array}
  \qquad
  V (w) = w \,,\qquad
  \rho_o (x)
  = \left(2+\sin (8\,x_2)\right) \caratt{[-1/4,1/4]^2} (x)
  \,,\qquad
  \begin{array}{r@{\,}c@{\,}l}
    r
    & =
    & 0.8\,,
    \\
    \ell
    & =
    & 1.5\,,
  \end{array}
\end{equation}
and as $\eta$ we choose~\eqref{eq:46}.

\begin{figure}[!h]
  \includegraphics[width=0.5\linewidth]{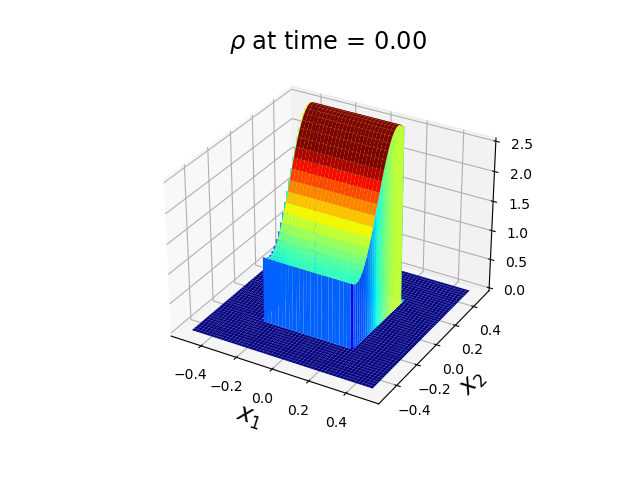}%
  \includegraphics[width=0.5\linewidth]{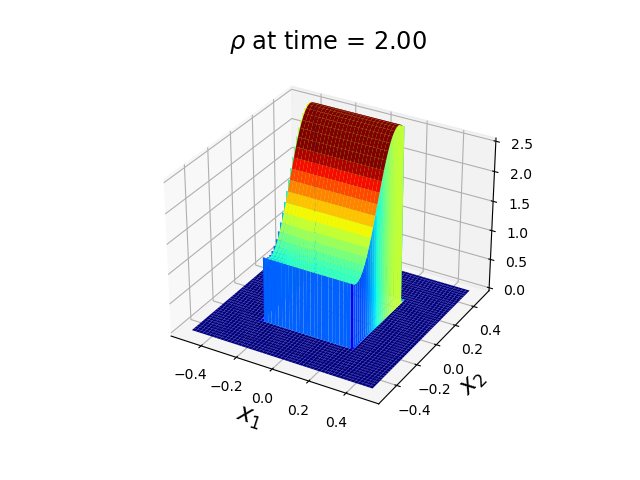}
  \caption{Numerical integration of equation~\eqref{eq:1} with
    parameters and initial datum as
    in~\eqref{eq:46}--\eqref{eq:45}. Coherently with
    Proposition~\ref{prop:stationary}, the resulting solution is
    stationary.}\label{fig:stationary}
\end{figure}

Note that $\eta$ is constant for $\norma{x} \leq r = 0.8$ and the
support of the initial datum has diameter $1/\sqrt2$. Since
$1/\sqrt2 < 0.8$, the resulting numerical integration yields a
stationary solution, coherently with
Proposition~\ref{prop:stationary}, see Figure~\ref{fig:stationary}.

Here, the numerical domain is $[-0.5, 0.5]^2$, the mesh consists of
$10^3 \times 10^3$ points and the upwind method was used.

\subsection{Encryption -- Decryption}
\label{subsec:encrypt-decrypt-imag}

Below, we write explicitly the dependence of the group $\mathcal{G}$
contructed in Theorem~\ref{thm:cauchyProblem} on $\eta$ and $V$ as
$\mathcal{G}^{\eta,V}$.

The well posedness proved in Theorem~\ref{thm:cauchyProblem}, and in
particular the reversibility shown therein, allows to consider the
following procedure:
\begin{enumerate}[label=\bf\arabic*.]
\item \label{item:28} Interpret $\rho_o = \rho_o (x)$ as information
  to be encrypted.
\item \label{item:27} Fix arbitrary $V$, $\eta$
  satisfying~\ref{item:7} and~\ref{item:6}. They are the encryption
  keys.
\item \label{item:29} Fix a positive $T$. Then, the encrypted signal
  is $\rho = \mathcal{G}^{\eta,V}_{\strut T} \rho_o$.\vspace{-12pt}
\item \label{item:8} To decrypt, integrate backwards since
  $\rho_o = \mathcal{G}^{\eta,V}_{\strut -T} \rho$.\vspace{-8pt}
\end{enumerate}

\noindent A result like the following one would play a key role in
assessing the reliability of the above encryption procedure.

\begin{problem}
  \label{pb:encrypt-decrypt}
  Call $\mathcal{K}_\eta$ and $\mathcal{K}_V$ the set of the keys
  $\eta$ and of $V$ that satisfy~\ref{item:6} and~\ref{item:7}. Find a
  sufficiently large class
  $\mathcal{R} \subset (\L1 \cap \L\infty \cap \BV) (\reali^n;
  \reali^n)$ such that for all $\rho_o \in \mathcal{R}$ the map
  $(\eta,V) \mapsto \mathcal{G}^{\eta,V}_{\strut T} \rho_o$ is
  surjective onto $\mathcal{R}$.
\end{problem}

Indeed, this surjectivity ensures that the above encryption procedure
is essentially unbreakable in the class $\mathcal{R}$. On the basis of
the properties of characteristic, it can be proved~\cite{KeimerPflug}
that $\mathcal{R}$ can not be chosen as large as
$(\L1 \cap \L\infty \cap \BV) (\reali^n; \reali^n)$. The introduction
of source terms~\cite{andrea1} is likely to allow for a larger class
$\mathcal{R}$.

At the numerical level, the above procedure needs an efficient
numerical algorithm to compute approximate solutions to~\eqref{eq:1}
that is \emph{reversible}, thus respecting this key property
of~\eqref{eq:1}.

\subsubsection{The Scalar 1-Dimensional Case}
\label{subsubsec:scalar-1-dimensional}

\begin{figure}[!h]
  \centering%
  \includegraphics[width=0.33\textwidth, trim=20 10 20 10
  keepaspectratio]{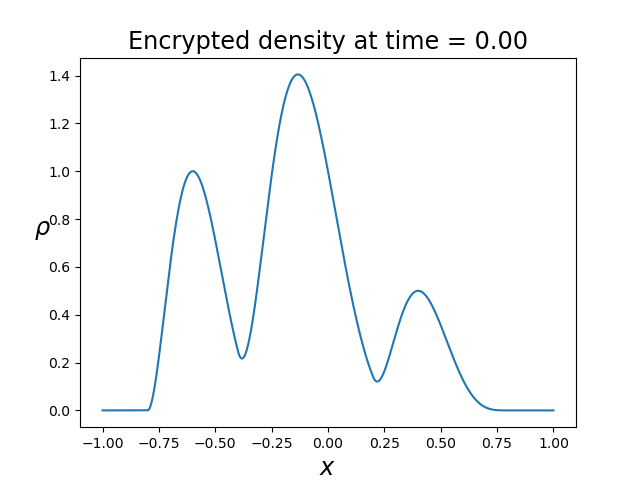}%
  \includegraphics[width=0.33\textwidth, trim=20 10 20 10
  keepaspectratio]{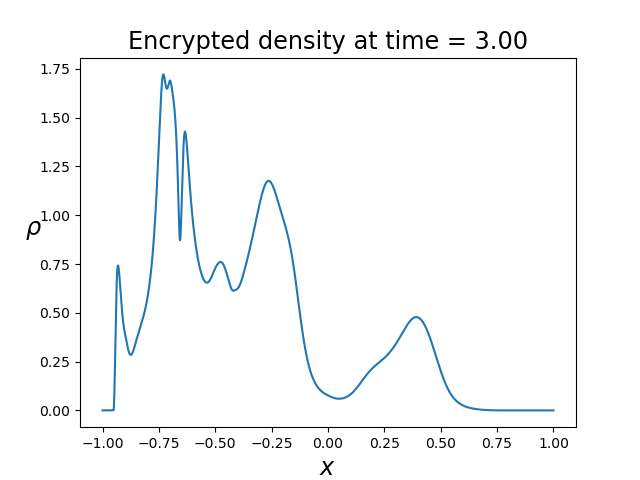}%
  \includegraphics[width=0.33\textwidth, trim=20 10 20 10
  keepaspectratio]{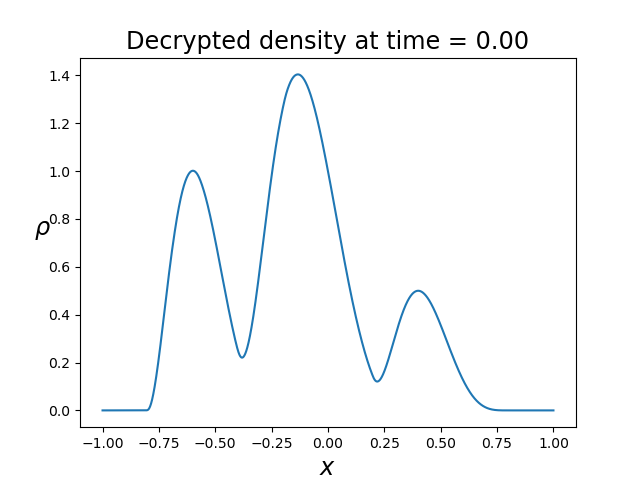}%
  \caption{Left, the original data~\eqref{eq:42}; center the encrypted
    data and, right, the decrypted data. Encryption is obtained
    through~\eqref{eq:1}--\eqref{eq:36}, decryption through an
    integration backward in time. Theorem~\ref{thm:cauchyProblem}
    guarantees the feasibility of this procedure.}
  \label{fig:crypt}
\end{figure}

Consider the encryption/decryption of, say, a signal. Indeed,
consider~\eqref{eq:1} in the simplest case where, with reference
to~\eqref{eq:1}, we set $n=1$, $m=1$ and
\begin{equation}
  \label{eq:36}
  \begin{array}{@{}r@{\,}c@{\,}l}
    n
    & =
    & 1
    \\
    m
    & =
    & 1
  \end{array}
  \qquad
  V (w) = \dfrac{w}{\sqrt{1+w^2}}
  \qquad
  \eta (\xi)
  =
  \xi \; \exp \left(- \dfrac{1}{1 - (\xi/\ell)^2}\right)
  \; \caratt{[-\ell, \ell]} (\xi)
  \qquad
  \ell = 1/4 \,.
\end{equation}
As initial datum, we choose
\begin{equation}
  \label{eq:42}
  \begin{array}{c}
    \rho_o (x) =
    \theta (x;\,-0.8,\,-0.2) +
    \theta (x;\,-0.4,\,0.4) +
    \theta (x;\, 0.2,\, 0.8)
    \\
    \mbox{ where } \quad
    \theta (x;\, a,\, b)
    =
    \left(1-\frac{x}{a}\right)^2 \; \left(1-\frac{x}{b}\right)^4 \;
    \caratt{[a,b]} (x)\,.
  \end{array}
\end{equation}
Figure~\ref{fig:crypt}
shows the initial datum $\rho_o$, corresponding to the original
signal, then the encrypted signal $\mathcal{G}_3 \rho_o$, i.e., the
solution to~\eqref{eq:1}--\eqref{eq:36}--\eqref{eq:42}, and the image
resulting from the decryption, i.e., the numerical approximation
$\tilde \rho$ of $\mathcal{G}_{-3} \, \mathcal{G}_3 \rho_o$. As
described in the proof of Theorem~\ref{thm:cauchyProblem}, decryption
is obtained through~\eqref{eq:1} where $V$ is replaced by $-V$.

A measure of the reliability of the encryption -- decryption process
is given by the $\L1$ difference between the original datum $\rho_o$
in~\eqref{eq:42} and $\tilde\rho$, see the values in~\eqref{eq:44}.
\begin{equation}
  \label{eq:44}
  \!\!\!
  \begin{array}{@{}r@{\,}c@{\,}l}
    \norma{\rho_o}_{\L1 (\reali; \reali)}
    & \approx
    & 0.8781547619047426
    \\
    \norma{\mathcal{G}_3\rho_o}_{\L1 (\reali; \reali)}
    & \approx
    & 0.8781547619047422
    \\
    \norma{\tilde\rho}_{\L1 (\reali; \reali)}
    & \approx
    & 0.8781547619047416
  \end{array}
  \quad
  \begin{array}{r@{\,}c@{\,}l@{}}
    \norma{\rho_o - \tilde\rho}_{\L1 (\reali; \reali)}
    & \approx
    & 8.880075473908466 \times 10^{-4}
    \\
    \dfrac{\norma{\rho_o - \tilde\rho}_{\L1 (\reali; \reali)}}{\norma{\rho_o}_{\L1 (\reali; \reali)}}
    & \approx
    & 0,00101121987366
  \end{array}
\end{equation}
This integration was performed on $[-1,1]$ with $10^5$ mesh points
using the upwind method.

As a side remark, note that the middle diagram in
Figure~\ref{fig:crypt} shows no clear formation of clusters,
coherently with the fact that $\eta$ changes sign, resulting in partly
attractive and partly repulsive interactions.

\subsubsection{$2$-Dimensonal Examples}

For completeness, we present two forward -- backward integrations
of~\eqref{eq:1} in the case $n=2$, $m=1$. We use the upwind scheme
with dimensional splitting with a regular $8000\times 8000$ mesh. The
non reversibility of the numerical method imposes to consider a small
time interval.

\begin{figure}[!h]
  \centering%
  \includegraphics[width=0.33\textwidth, trim=20 10 20 10
  keepaspectratio]{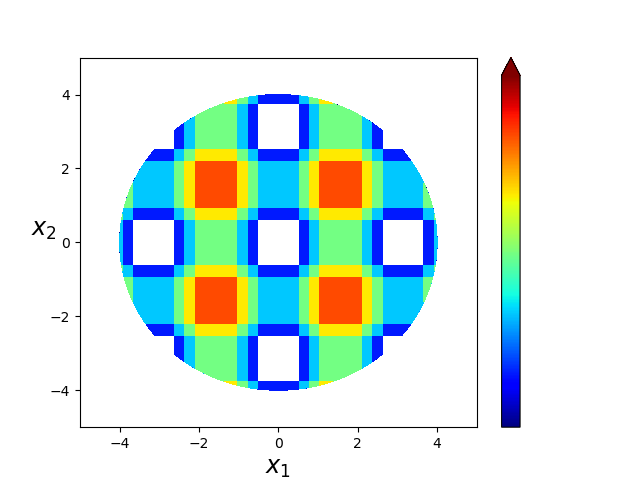}%
  \includegraphics[width=0.33\textwidth, trim=20 10 20 10
  keepaspectratio]{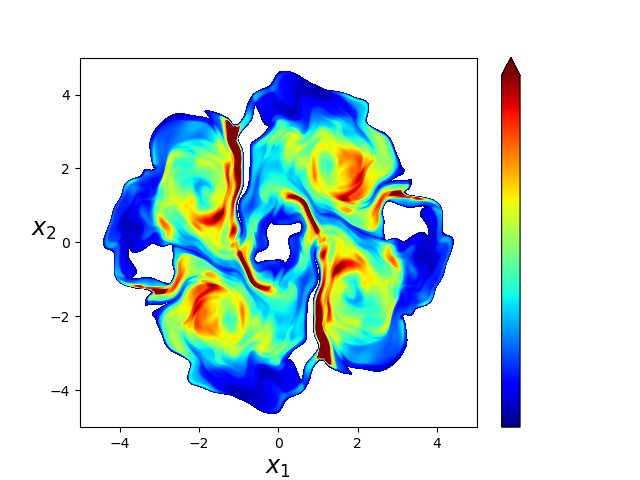}%
  \includegraphics[width=0.33\textwidth, trim=20 10 20 10
  keepaspectratio]{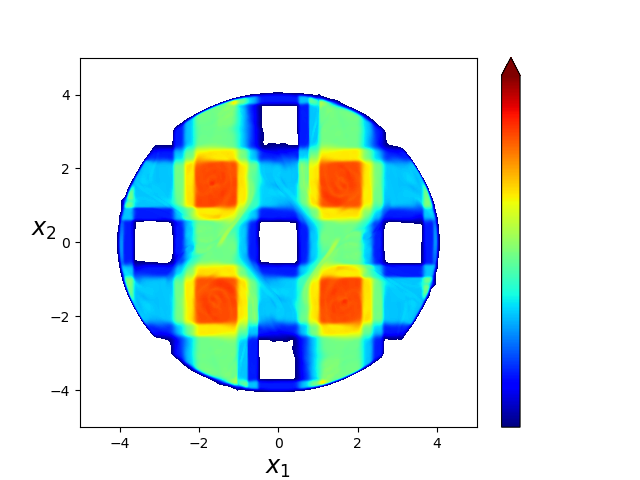}%
  \caption{From left to right, the original $\rho_o$, its encryption
    and the decryption resulting from backward integration
    of~\eqref{eq:1}--\eqref{eq:10}--\eqref{eq:30}.}
  \label{fig:simple03}
\end{figure}
\begin{figure}[!h]
  \centering%
  \includegraphics[width=0.33\textwidth, trim=20 10 20 10
  keepaspectratio]{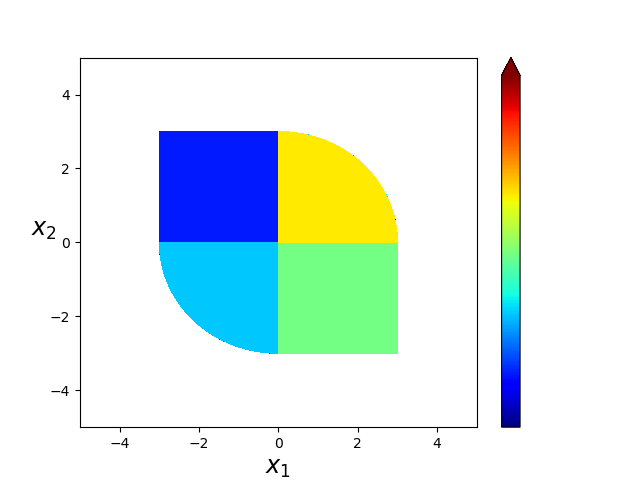}%
  \includegraphics[width=0.33\textwidth, trim=20 10 20 10
  keepaspectratio]{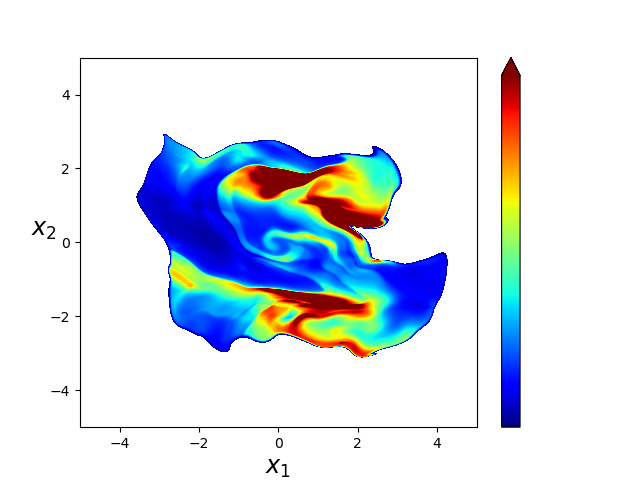}%
  \includegraphics[width=0.33\textwidth, trim=20 10 20 10
  keepaspectratio]{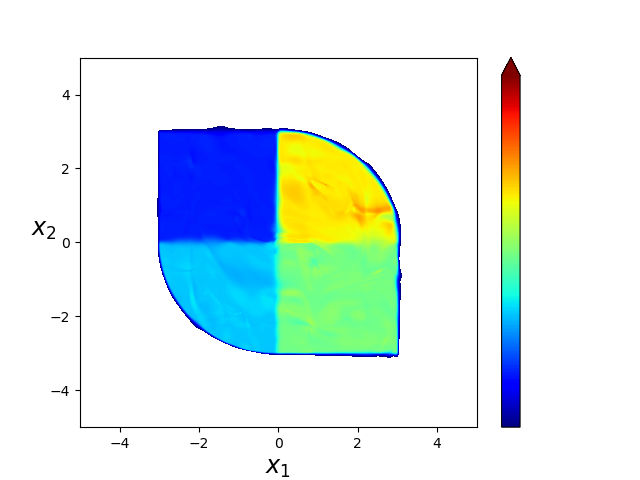}%
  \caption{From left to right, the original $\rho_o$, its encryption
    and the decryption resulting from backward integration
    of~\eqref{eq:1}--\eqref{eq:10}--\eqref{eq:37}.}
  \label{fig:simple02}
\end{figure}

In~\eqref{eq:1}, we set
\begin{equation}
  \label{eq:10}
  \begin{array}{@{}r@{\,}c@{\,}l}
    n
    & =
    & 2
    \\
    m
    & =
    & 1
    \\
    T
    & =
    & 0.25
  \end{array}
  \qquad
  V (w) = \left[{0 \;-1 \atop 1 \quad 0}\right] \dfrac{w}{\sqrt{1+\norma{w}^2}}
  \qquad
  \eta (\xi)
  =
  \cos \left(\frac\pi2 \; \norma{\xi}^2\right) \; \caratt{\norma{\xi}<1} (x)\,.
\end{equation}
and choose the initial data
\begin{eqnarray}
  \label{eq:30}
  \!\!\!\!\!\!\!\!\!\!\!\!\!\!\!\!\!
  \rho_o (x)
  & \!\!\!\!\!=\!\!\!\!\!
  & \left(
    \lfloor 4 \sin^2 x_1 \rfloor + \lfloor 3 \sin^2 x_2 \rfloor \right)
    \; \caratt{\{\norma{x}<4\}} (x) ;
  \\
  \label{eq:37}
  \!\!\!\!\!\!\!\!\!\!\!\!\!\!\!\!\!
  \rho_o (x)
  & \!\!\!\!\!=\!\!\!\!\!
  & 4 \caratt{\{x_1>0; x_2 >0; \norma{x}<3\}} (x) {+}
    \caratt{[-3,0] {\times} [0,3]} (x) {+}
    2 \caratt{\{x_1<0; x_2 <0; \norma{x}<3\}} (x) {+}
    3\caratt{[0,3] {\times} [-3,0]} (x)
\end{eqnarray}
where $\lfloor \,\cdot\,\rfloor$ is the integer part. The numerical solutions are in Figure~\ref{fig:simple03} and in Figure~\ref{fig:simple02}.


As a side remark, observe that a possibly
better encryption can be obtained extending~\eqref{eq:1} adding
suitable source terms, see~\cite{andrea1}.

\section{Proofs}
\label{sec:technical-proofs}


We set $\reali_+ = [0, +\infty\mathclose[$. For
$n \in \naturali \setminus \{0\}$, $\reali^n$ is equipped with the
Euclidean norm: this choice applies also to matrix spaces such as
$\reali^{n\times m}$. For $x \in \reali^n$ and $r>0$, $B (x,r)$ is the
open ball in $\reali^n$ centered at $x$ with radius $r$. Throughout,
$\nabla$ is the usual differential operator with respect to the
spatial coordinates, while $DV$ is the total derivative of the map $V$
with respect to all its variables.

$\mathcal{L}$ is the Lebesgue measure in $\reali^n$. Concerning
function spaces and norms, we convene that for a measurable function
$\rho \colon \reali^n \to \reali^m$ and for any $p \in \{1, +\infty\}$
\begin{displaymath}
  \begin{array}{@{}l@{}}
    \norma{\rho}_{\L1 (\reali^n; \reali^m)}
    =
    \int_{\reali^n} \norma{\rho (x)} \d{x}
    \\
    \norma{\rho}_{\L\infty (\reali^n; \reali^m)}
    =
    \esssup_{x \in \reali^n} \norma{\rho (x)}
  \end{array}
  \;\mbox{ and }\;
  \begin{array}{@{}l@{}}
    \norma{\rho}_{\W1p (\reali^n; \reali^m)}
    =
    \norma{\rho}_{\L{p} (\reali^n; \reali^m)}
    +
    \norma{D \rho}_{\L{p} (\reali^n; \reali^{n\times m})}
    \\
    \norma{\rho}_{\C0 (\reali_+; \L{p} (\reali^n;\reali^m))}
    =
    \sup_{t \in \reali_+} \norma{\rho (t)}_{\L{p} (\reali^n;\reali^m)} \,.
  \end{array}
\end{displaymath}
By $\BV (\reali^n; \reali)$ we mean the set of functions in
$\Lloc1(\reali^n; \reali)$ such that their distributional derivative
is a (vector) Radon measure with finite total variation. In
particular, a function in $\BV (\reali^n; \reali)$ needs not be also
in $\L1 (\reali^n; \reali)$, as in the case of the function
$\rho \equiv 1$ on $\reali^n$.

In connection with~\eqref{eq:23}, for $i \in \{1, \ldots, m\}$,
introduce the characteristics
$X_i \colon \reali \times \reali \times \reali^n \to \reali^n$ where
\begin{equation}
  \label{eq:29}
  t \mapsto X_i (t;t_o,x_o)
  \quad \mbox{ solves } \quad
  \left\{
    \begin{array}{@{}l@{}}
      \dot x = v_i (t,x)
      \\
      x (t_o) = x_o \,.
    \end{array}
  \right.
\end{equation}
It is immediate to check that when $v_i$ is as in~\eqref{eq:25}, for
any $\rho \in \C0\left(\reali;\L1 (\reali^n; \reali^m)\right)$,
assumptions~\ref{item:6} and~\ref{item:7} ensure the global existence
and regularity of $X_i$.

For any $\rho_o \in \L1 (\reali^n; \reali^m)$, recall that the
solution to~\eqref{eq:23} with initial datum $\rho_o$ assigned at time
$t_o=0$ is given, for $i \in \{1, \ldots, m\}$, through the Lagrangian
representation
\begin{equation}
  \label{eq:7}
  \rho_i (t,x)
  =
  \rho_{o,i} \left(X_i (0;t,x)\right)
  \exp \left(
    -\int_0^t (\nabla \cdot v_i) \left(\tau; X_i (\tau;t,x)\right) \d\tau
  \right)
\end{equation}
which is a distributional solution and also a {Kru\v zkov} solution in
the sense of~\cite[Definition~1]{MR0267257}. Indeed, several results
in the literature, see for instance~\cite[Lemma~5]{MR4371486} or
\cite[Corollary II.1]{zbMATH04140224}, ensure that in the case of the
Cauchy problem for~\eqref{eq:23}, the concepts of \emph{weak} and
\emph{entropy} (or \emph{Kru\v zkov}) solutions coincide. Under
assumptions~\ref{item:6} and~\ref{item:7}, when $\rho_{o,i}$ is
sufficiently regular, \eqref{eq:7} also yields the \emph{strong}
solution to~\eqref{eq:23} with $v$ as in~\eqref{eq:25}.

For $T > 0$, introduce the set
\begin{equation}
  \label{eq:6}
  \mathcal{V}_T
  =
  \left\{
    v \in (\C0 \cap \L\infty) \left([0,T]\times\reali^n; \reali^{n\times m}\right)
    \colon
    \begin{array}{@{}l@{}}
      v (t)
      \in
      \C2 (\reali^n; \reali^{n\times m})
      \quad     \forall\, t \in [0,T]
      \\
      \nabla v
      \in
      \L\infty ([0,T] \times \reali^n; \reali^{n^2\times m})
      \\
      \nabla (\nabla \cdot v)
      \in
      \L1 ([0,T] \times \reali^n; \reali^{n\times m})
    \end{array}
  \right\}
\end{equation}
and, for a fixed $\rho_o \in \L1 (\reali^n; \reali^m)$, the map
\begin{equation}
  \label{eq:2}
  \begin{array}[b]{@{}c@{\,}c@{\,}c@{\,}c@{\,}c@{}}
    \Sigma
    & \colon
    & \mathcal{V}_T
    & \to
    & \C0\left([0,T]; \L1 (\reali^n; \reali^m)\right)
    \\
    &
    & v
    & \mapsto
    & \rho
  \end{array}
  \quad
  \mbox{where }
  \quad
  \left\{
    \begin{array}{l}
      \partial_t \rho_i + \nabla \cdot ( \rho_i \; v_i) = 0
      \quad i\in \{1, \ldots,m\}
      \\
      \rho (0) = \rho_o \,.
    \end{array}
  \right.
\end{equation}
\noindent The Cauchy problem in~\eqref{eq:2} consists of $m$ linear,
scalar, independent conservation laws. Thus, the proof of its well
posedness is very similar to that of various results in the
literature. Here we refer to~\cite{MR4371486} and provide precise
references.

\begin{lemma}
  \label{lem:solution}
  Fix $T>0$ and let $\rho_o \in \L1 (\reali^n; \reali^m)$. Then, the
  map $\Sigma$ in~\eqref{eq:2} is well defined and enjoys the
  following properties:
  \begin{enumerate}[label=($\mathbf{\Sigma}\textbf{\arabic*}$)]
  \item \label{item:1} For all $t \in [0,T]$:
    \begin{equation}
      \label{eq:13}
      \norma{(\Sigma v) (t)}_{\L1 (\reali^n; \reali^m)}
      =
      \norma{\rho_o}_{\L1 (\reali^n;\reali^m)} \,.
    \end{equation}

  \item \label{item:15} If $\rho_o \in \L\infty (\reali^n; \reali^m)$,
    for all $t \in [0,T]$:
    \begin{equation}
      \label{eq:14}
      \norma{(\Sigma v) (t)}_{\L\infty (\reali^n; \reali^m)}
      \leq
      \norma{\rho_o}_{\L\infty (\reali^n;\reali^m)} \;
      \exp\left(
        \norma{\nabla \cdot v}_{\L1 ([0,t];\L\infty (\reali^n;\reali^m))}
      \right) \,.
    \end{equation}

  \item \label{item:16} If $\rho_o \in \BV (\reali^n; \reali^m)$, for
    all $t \in [0,T]$:
    \begin{eqnarray}
      \label{eq:15}
      \tv\left((\Sigma v) (t)\right)
      & \leq
      & \left(
        1
        +
        \norma{\nabla (\nabla \cdot v)}_{\L1 ([0,t];\L1 (\reali^n; \reali^{n\times m}))}
        \right)
      \\
      \nonumber
      &
      & \times
        \exp
        \norma{\nabla v}_{\L1 ([0,t];\L\infty (\reali^n; \reali^{n^2\times m}))}
        \left(
        \norma{\rho_o}_{\L\infty (\reali^n; \reali^m)} + \tv (\rho_o)
        \right) \,.
    \end{eqnarray}

  \item \label{item:2} If
    $\rho_o \in (\L\infty \cap \BV) (\reali^n; \reali^m)$,
    $v_1,v_2 \in \mathcal{V}_T$ and
    $\nabla \cdot (v_2 - v_1) \in {\L1 ([0,t]\times \reali^n;
      \reali^m)}$ for a $t \in [0,T]$, then
    \begin{eqnarray*}
      &
      &
        \norma{(\Sigma v_2 - \Sigma v_1) (t)}_{\L1 (\reali^n; \reali^m)}
      \\
      & \leq
      & \exp \int_0^t
        \max\left\{
        \norma{\nabla v_1 (\tau)}_{\L\infty (\reali^n; \reali^{n^2\times m})}
        ,
        \norma{\nabla v_2 (\tau)}_{\L\infty (\reali^n; \reali^{n^2\times m})}
        \right\} \d\tau
      \\
      &
      & \quad\times
        \left(
        1
        +
        \int_0^t
        \max\left\{
        \norma{\nabla \left(\nabla \cdot v_1 (\tau)\right)}_{\L1 (\reali^n; \reali^{n\times m})}
        ,
        \norma{\nabla \left(\nabla \cdot v_2 (\tau)\right)}_{\L1 (\reali^n; \reali^{n\times m})}
        \right\}\d\tau
        \right)
      \\
      &
      & \quad\times
        \left(\norma{\rho_o}_{\L\infty (\reali^n; \reali^m)} + \tv (\rho_o)\right) \;
        \norma{v_2 - v_1}_{\L1([0,t];\L\infty(\reali^n; \reali^{n\times m})}
      \\
      &
      & + \norma{\rho_o}_{\L\infty (\reali^n; \reali^m)} \;
        \norma{\nabla \cdot v_2 - \nabla\cdot v_1}_{\L1 ([0,t]\times \reali^n; \reali^m)} \,.
    \end{eqnarray*}

  \item \label{item:3} If
    $\rho_o \in (\L\infty \cap \BV) (\reali^n; \reali^m)$, for all
    $t_1,t_2 \in [0,T]$,
    \begin{eqnarray*}
      \norma{(\Sigma v) (t_2) - (\Sigma v) (t_1)}_{\L1 (\reali^n; \reali^m)}
      & \leq
      & \left(
        1
        +
        \norma{\nabla (\nabla \cdot v)}_{\L1 ([0,t_1 \vee t_2];\L1 (\reali^n; \reali^{n\times m}))}
        \right)
      \\
      &
      & \times
        \exp
        \norma{\nabla v}_{\L1 ([0,t_1 \vee t_2];\L\infty (\reali^n; \reali^{n^2\times m}))}
      \\
      &
      & \times
        \left(
        \norma{\rho_o}_{\L\infty (\reali^n; \reali^m)}
        +
        \tv (\rho_o)
        \right)
        \modulo{t_2 - t_1} \,.
    \end{eqnarray*}

  \item \label{item:5} If $\rho_{o,i} \geq 0$ for an
    $i \in \{1, \ldots,n\}$, then for all $t \in [0,T]$,
    $\left((\Sigma v) (t)\right)_i \geq 0$.

  \item \label{item:13} Let
    $\widehat \rho_o, \widecheck \rho_o \in \L1 (\reali^n; \reali^m)$
    and call $\widehat\Sigma, \widecheck\Sigma$ the corresponding maps
    defined as in~\eqref{eq:2} with
    $\widehat \rho_o, \widecheck \rho_o$ as initial datum. Then, for
    all $t \in [0,T]$,
    \begin{displaymath}
      \norma{(\widehat\Sigma v) (t)
        - (\widecheck\Sigma v) (t)}_{\L1 (\reali^n; \reali^m)}
      =
      \norma{\widehat \rho_o - \widecheck \rho_o}_{\L1 (\reali^n;\reali^m)} \,
    \end{displaymath}
    and if moreover
    $\widehat \rho_o, \widecheck \rho_o \in \L\infty (\reali^n;
    \reali^m)$, for all $t \in [0,T]$,
    \begin{displaymath}
      \norma{(\widehat\Sigma v) (t)
        - (\widecheck\Sigma v) (t)}_{\L\infty (\reali^n; \reali^m)}
      \leq
      \norma{\widehat \rho_o - \widecheck \rho_o}_{\L\infty (\reali^n;\reali^m)}
      \,
      \exp
      \left(
        \norma{\nabla \cdot v}_{\L1 ([0,t];\L\infty (\reali^n;\reali^m))}
      \right) \,.
    \end{displaymath}
  \end{enumerate}
\end{lemma}

\begin{proofof}{Lemma~\ref{lem:solution}}
  The equality in~\ref{item:1} directly follows from~\eqref{eq:7},
  also using the change of coordinates $\xi = X_i (0;t,x)$ whose
  Jacobian determinant $J_i = J_i (t,x)$ satisfies
  $J_i (t,x) = \exp \int_0^t \nabla \cdot v_i\left(\tau, X_i
    (\tau;0,x)\right) \d\tau$, see~\cite[(H3) in
  Proposition~3]{MR4371486} for the details. The bound
  in~\ref{item:15} follows from a direct estimate based
  on~\eqref{eq:7}, see~\cite[(H4) in Proposition~3]{MR4371486}. The
  bound~\ref{item:16} is proved through an approximation argument,
  see~\cite[(H6) in Proposition~3]{MR4371486} for the details.  The
  estimate in~\ref{item:2} follows from~\cite[(H5) in
  Proposition~3]{MR4371486}. The continuity in time~\ref{item:3}
  is~\cite[(H7) in Proposition~3]{MR4371486}. The proof
  of~\ref{item:5} is immediate, thanks to~\eqref{eq:7}. Finally,
  \ref{item:13} directly follows from~\ref{item:1}, thanks to the
  linearity of the differential equation in~\eqref{eq:2}.
\end{proofof}

For a $V$ as in~\ref{item:7}, define the maps
\begin{equation}
  \label{eq:5}
  \begin{array}[b]{@{}c@{\,}c@{\,}c@{\,}c@{\,}c@{}}
    \tilde\Pi
    & \colon
    & \L1 (\reali^n; \reali^m)
    & \to
    & \C2 (\reali^n; \reali^{n\times m})
    \\
    &
    & \rho
    & \mapsto
    & V\left(\nabla (\rho * \eta)\right)
  \end{array}
  \mbox{ where }
  \left(\left(\nabla (\rho *\eta)\right) \!(x)\right)_{ji}
  =
  \int_{\reali^n} \rho_i (x-\xi)\; \partial_{x_j}\eta (\xi) \, \d\xi
\end{equation}
and for a $T>0$
\begin{equation}
  \label{eq:47}
  \begin{array}[b]{ccccc}
    \Pi
    & \colon
    & \C0\left([0,T];\L1 (\reali^n;\reali^m) \right)
    & \to
    & \L\infty \left([0,T]; \L\infty(\reali^n; \reali^m)\right)
    \\
    &
    & \rho
    & \mapsto
    & \left[ t \mapsto \tilde\Pi\left(\rho (t)\right) \right] \,.
  \end{array}
\end{equation}

\begin{lemma}
  \label{lem:product}
  Fix $p \in [1, +\infty]$, $V$ satisfying~\ref{item:7} and let
  \begin{equation}
    \label{eq:40}
    \eta
    \in
    (\C3 \cap \W3{p}) (\reali^n; \reali)
    \quad \mbox{ and } \quad
    \nabla \eta \in \L\infty (\reali^n; \reali^n)
    \quad \mbox{ and } \quad
    \Delta \eta \in \L\infty (\reali^n; \reali) \,.
  \end{equation}
  Then, $\tilde\Pi$ in~\eqref{eq:5} is well defined and for all
  $\rho \in \L1 (\reali^n; \reali^m)$
  \begin{eqnarray}
    \label{eq:3}
    \tilde\Pi\rho
    & \in
    & \C2 (\reali^n; \reali^{n\times m})
    \\
    \label{eq:16}
    \norma{\tilde\Pi\rho}_{\L{p} (\reali^n; \reali^{n\times m})}
    & \leq
    & L_V \; \norma{\nabla \eta}_{\L{p} (\reali^n; \reali^n)} \;
      \norma{\rho}_{\L1 (\reali^n;\reali^m)}
    \\
    \label{eq:17}
    \norma{\nabla \cdot \tilde\Pi\rho}_{\L{p} (\reali^n; \reali^m)}
    & \leq
    & L_V \; \norma{\Delta \eta}_{\L{p} (\reali^n; \reali)} \;
      \norma{\rho}_{\L1 (\reali^n;\reali^m)}
    \\
    \label{eq:18}
    \!\!\!\!
    \norma{\nabla \tilde\Pi\rho}_{\L{p} (\reali^n; \reali^{n^2\times m})}
    & \leq
    & L_V \; \norma{D^2\eta}_{\L{p} (\reali^n; \reali^{n\times n})} \;
      \norma{\rho}_{\L1 (\reali^n;\reali^m)}
    \\
    \label{eq:19}
    \norma{\nabla(\nabla \cdot \tilde\Pi\rho)}_{\L{p} (\reali^n; \reali^{n\times m})}
    & \leq
    & L_V \,
      \norma{D^2 \eta}_{\L{p} (\reali^n; \reali^{n\times n})} \,
      \norma{\Delta\eta}_{\L\infty (\reali^n; \reali)} \,
      \norma{\rho}_{\L1 (\reali^n;\reali^m)}^2
    \\
    \nonumber
    &
    & +
      L_V \,
      \norma{D^3 \eta}_{\L{p} (\reali^n; \reali^{n\times n\times n})} \,
      \norma{\rho}_{\L1 (\reali^n;\reali^m)}
      .
  \end{eqnarray}
  Moreover, for all $\rho_1,\rho_2 \in \L1 (\reali^n; \reali^m)$,
  setting
  \begin{displaymath}
    C_\Pi = L_V \left(1 + \min \left\{\norma{\rho_1}_{\L1 (\reali^n;
          \reali^m)}, \norma{\rho_2}_{\L1 (\reali^n; \reali^m)}\right\}
      \norma{\nabla \eta}_{\L{p} (\reali^n; \reali^n)}\right)
  \end{displaymath}
  we have
  \begin{eqnarray}
    \label{eq:20}
    \norma{\tilde\Pi (\rho_2) - \tilde\Pi (\rho_1)}_{{\L{p} (\reali^n; \reali^{n\times m})}}
    & \leq
    & L_V \; \norma{\nabla \eta}_{\L{p} (\reali^n; \reali^n)} \;
      \norma{\rho_2 - \rho_1}_{\L1 (\reali^n;\reali^m)}
    \\
    \label{eq:22}
    \norma{\nabla \cdot \tilde\Pi (\rho_2) - \nabla \cdot \tilde\Pi (\rho_1)}_{{\L{p} (\reali^n; \reali^m)}}
    & \leq
    & C_\Pi \; \norma{\Delta \eta}_{\L{p} (\reali^n; \reali)} \;
      \norma{\rho_2 - \rho_1}_{\L1 (\reali^n;\reali^m)}
    \\
    \label{eq:11}
    \norma{\nabla \tilde\Pi (\rho_2) - \nabla \tilde\Pi (\rho_1)}_{{\L{p} (\reali^n; \reali^{n^2\times m})}}
    & \leq
    & C_\Pi \; \norma{D^2\eta}_{\L{p} (\reali^n; \reali^{n\times n})} \;
      \norma{\rho_2 - \rho_1}_{\L1 (\reali^n;\reali^m)} \,.
  \end{eqnarray}
\end{lemma}

Observe that Lemma~\ref{lem:product} admits a straightforward
extension to $\Pi$ as defined in~\eqref{eq:47}, where
$\norma{\;\cdot\;}_{\L1 (\reali^n;\reali^m)}$ is replaced by
$\norma{\;\cdot\;}_{\C0([0,T];\L1 (\reali^n;\reali^m))}$,
see~\cite{Claudia01}.

\begin{proofof}{Lemma~\ref{lem:product}}
  The regularity~\eqref{eq:3} and all subsequent estimates are
  immediate consequences of~\ref{item:7} and of the classical
  properties of the convolution product, see for
  instance~\cite[\S~IV.4]{Brezis}.
\end{proofof}

\begin{proofof}{Theorem~\ref{thm:cauchyProblem}}
  Fix a positive $T$ and a
  $\rho_o \in (\L1 \cap \L\infty \cap \BV) (\reali^n;
  \reali^m)$. Introduce the set
  \begin{equation}
    \label{eq:4}
    \mathcal{R}
    =
    \left\{
      \rho \in \C0 \left([0,T]; \L1 (\reali^n; \reali^m)\right)
      \colon
      \begin{array}{@{}l@{}}
        \rho (0) = \rho_o \mbox{ and }
        \\
        \forall \, t \in [0,T] \quad
        \norma{\rho (t)}_{\L1(\reali^n;\reali^m)} = \norma{\rho_o}_{\L1(\reali^n;\reali^m)}
      \end{array}
    \right\}
  \end{equation}
  and define the composition
  \begin{equation}
    \label{eq:8}
    \begin{array}{ccccc}
      \mathcal{T}
      & \colon
      & \mathcal{R}
      & \to
      & \mathcal{R}
      \\
      &
      & \rho
      & \mapsto
      & (\Sigma \circ \Pi) (\rho)
    \end{array}
  \end{equation}
  with $\Sigma$ defined in~\eqref{eq:2} and $\Pi$ in~\eqref{eq:47}. We
  now prove that $\mathcal{T}$ is well defined and a
  contraction. Indeed, by Definition~\ref{def:solution} a map $\rho$
  is a solution to~\eqref{eq:1} if and only if it is a fixed point of
  $\mathcal{T}$ on the time interval where it is defined.

  \paragraph{Claim~1: $\mathcal{T}$ is well defined.}

  Prove first that if $\rho \in \mathcal{R}$ then
  $\Pi\rho \in \mathcal{V}_T$, with $\mathcal{V}_T$ as
  in~\eqref{eq:6}. For any $\rho \in \mathcal{R}$,
  $\Pi\rho \in \L\infty ([0,T]\times\reali^n;\reali^{n\times m})$
  by~\eqref{eq:16}.  To prove the continuity of $\Pi \rho$, fix
  $\left(\bar t, \bar x\right) \in [0, T] \times \reali^n$ and
  $\eps > 0$. Then, since $\rho \in \mathcal R$ and~\eqref{eq:3},
  there exists $\delta > 0$ such that
  \begin{equation*}
    \norma{\rho(t) - \rho(\bar t)}_{\LL1\left(\reali^n; \reali^m\right)} <
    \frac{\eps}{2 L_V \norma{\nabla \eta}_{\LL\infty\left(\reali^n; \reali\right)}}
  \end{equation*}
  for every $t \in [0, T]$, $\modulo{t - \bar t} < \delta$ and
  \begin{equation*}
    \norma{\tilde \Pi\left(\rho(\bar t)\right)(x)
      - \tilde \Pi\left(\rho(\bar t)\right)(\bar x)}
    < \frac{\eps}{2}
  \end{equation*}
  for every $\norma{x - \bar x} < \delta$. Hence, using~\eqref{eq:20},
  we deduce that
  \begin{eqnarray*}
    \norma{ \Pi \rho(t)(x) - \Pi \rho(\bar t) (\bar x)}
    & \le
    &
      \norma{ \Pi \rho(t)(x) - \Pi \rho(\bar t) (x)} +
      \norma{ \Pi \rho(\bar t)(x) - \Pi \rho(\bar t) (\bar x)}
    \\
    & \le
    &
      \norma{ \tilde \Pi \left(\rho(t)\right)(x) -
      \tilde \Pi \left(\rho(\bar t)\right) (x)}
      + \norma{ \tilde \Pi \left(\rho(\bar t)\right) (x)
      - \tilde \Pi \left(\rho(\bar t)\right) (\bar x)}
    \\
    & \le
    &
      \norma{ \tilde \Pi \left(\rho(t)\right) -
      \tilde \Pi \left(\rho(\bar t)\right)}_{\LL\infty\left(\reali^n;
      \reali^{n\times m}\right)}
      + \frac{\eps}{2}
    \\
    & \le
    &
      L_V \norma{\nabla \eta}_{\LL\infty\left(\reali^n; \reali^n\right)}
      \norma{ \rho(t) - \rho(\bar t)}_{\LL1\left(\reali^n;
      \reali^{m}\right)}
      + \frac{\eps}{2} < \eps,
  \end{eqnarray*}
  proving that
  $\Pi \rho \in \C0\left([0, T] \times \reali^n; \reali^{m \times
      n}\right)$.
  The requirement
  $\nabla \Pi\rho \in \L\infty ([0,T] \times \reali^n; \reali^{m\times
    n\times n})$ follows from~\eqref{eq:18} and~\eqref{eq:4}, while
  $\nabla (\nabla \cdot \Pi\rho) \in \L1([0,T] \times \reali^n;
  \reali^{n\times m})$ by~\eqref{eq:19}.

  Prove now that if $v \in \mathcal{V}_T$, then
  $\Sigma v \in \mathcal{R}$, with $\mathcal{R}$ as
  in~\eqref{eq:4}. Simply apply~\eqref{eq:13} in~\ref{item:1}
  and~\ref{item:3} using
  $\rho_o \in (\L1 \cap \L\infty \cap \BV) (\reali^n; \reali^m)$.

  Claim~1 is proved.

  \paragraph{Claim~2: For $T$ small, $\mathcal{T}$ is a contraction.}
  Fix $\rho_1,\rho_2$ in $\mathcal{R}$ and call $v_i = \Pi \rho_i$ for
  $i\in \{1,2\}$. Then, using~\eqref{eq:22}, for $t \in [0,T]$,
  \begin{eqnarray*}
    &
    & \norma{\nabla \cdot v_2 (t) - \nabla \cdot v_1 (t)}_{\L1 (\reali^n; \reali^m)}
    \\
    & \leq
    & L_V \left(1 +
      \norma{\rho_o}_{\L1 (\reali^n; \reali^m)} \;
      \norma{\nabla \eta}_{\L1 (\reali^n; \reali^n)}\right)
      \norma{\Delta \eta}_{\L1 (\reali^n; \reali)} \;
      \norma{\rho_2 (t) - \rho_1 (t)}_{\L1 (\reali^n;\reali^m)}
  \end{eqnarray*}
  which is finite, so that we can apply~\ref{item:2} in
  Lemma~\ref{lem:product} and obtain:
  \begin{eqnarray*}
    &
    & \norma{(\mathcal{T}\rho_2) (t) - (\mathcal{T}\rho_1) (t)}_{\L1 (\reali^n;\reali^m)}
    \\
    & =
    & \norma{(\Sigma v_2) (t) - (\Sigma v_1) (t)}_{\L1 (\reali^n;\reali^m)}
    \\
    & \leq
    & \exp \int_0^t
      \max\left\{
      \norma{\nabla v_1 (\tau)}_{\L\infty (\reali^n; \reali^{n^2\times m})}
      ,
      \norma{\nabla v_2 (\tau)}_{\L\infty (\reali^n; \reali^{n^2\times m})}
      \right\} \d\tau
    \\
    &
    & \quad\times
      \left(
      1
      +
      \int_0^t
      \max\left\{
      \norma{\nabla \left(\nabla \cdot v_1 (\tau)\right)}_{\L1 (\reali^n; \reali^{n\times m})}
      ,
      \norma{\nabla \left(\nabla \cdot v_2 (\tau)\right)}_{\L1 (\reali^n; \reali^{n\times m})}
      \right\}\d\tau
      \right)
    \\
    &
    & \quad\times
      \left(\norma{\rho_o}_{\L\infty (\reali^n; \reali^m)} + \tv (\rho_o)\right) \;
      \norma{v_2 - v_1}_{\L1([0,t];\L\infty(\reali^n; \reali^{n\times m})}
    \\
    &
    & + \norma{\rho_o}_{\L\infty (\reali^n; \reali^m)} \;
      \norma{\nabla \cdot v_2 - \nabla\cdot v_1}_{\L1 ([0,t] \times \reali^n; \reali^m)}
    \\
    & \leq
    & \exp \int_0^t
      \max\left\{
      \norma{\nabla \Pi\rho_1 (\tau)}_{\L\infty (\reali^n; \reali^{n^2\times m})}
      ,
      \norma{\nabla \Pi\rho_2 (\tau)}_{\L\infty (\reali^n; \reali^{n^2\times m})}
      \right\} \d\tau
    \\
    &
    & \quad\times
      \left(
      1
      +
      \int_0^t
      \max\left\{
      \norma{\nabla \left(\nabla \cdot \Pi\rho_1 (\tau)\right)}_{\L1 (\reali^n; \reali^{n\times m})}
      ,
      \norma{\nabla \left(\nabla \cdot \Pi\rho_2 (\tau)\right)}_{\L1 (\reali^n; \reali^{n\times m})}
      \right\}\d\tau
      \right)
    \\
    &
    & \quad\times
      \left(\norma{\rho_o}_{\L\infty (\reali^n; \reali^m)} + \tv (\rho_o)\right) \;
      \norma{\Pi\rho_2 - \Pi\rho_1}_{\L1([0,t];\L\infty(\reali^n; \reali^{n\times m})}
    \\
    &
    & + \norma{\rho_o}_{\L\infty (\reali^n; \reali^m)} \;
      \norma{\nabla \cdot \Pi\rho_2 - \nabla\cdot \Pi\rho_1}_{\L1 ([0,t] \times \reali^n; \reali^m)}
    \\
    & \leq
    & \exp
      \left(
      L_V \norma{D^2 \eta}_{\L\infty (\reali^n;\reali^{n\times n})}
      \int_0^t
      \max\left\{
      \norma{\rho_1 (\tau)}_{\L1 (\reali^n; \reali^m)}
      ,
      \norma{\rho_2 (\tau)}_{\L1 (\reali^n; \reali^m)}
      \right\} \d\tau
      \right)
    \\
    &
    & \quad\times
      \left(
      1
      +
      L_V
      \left(
      \norma{D^2 \eta}_{\L\infty (\reali^n;\reali^{n\times n})}
      \norma{\Delta\eta}_{\L\infty (\reali^n; \reali)}
      + \norma{D^3 \eta}_{\L\infty (\reali^n;\reali^{n\times n\times n})}
      \right)
      \right.
    \\
    &
    & \quad\times
      \left.
      \int_0^t
      \max\left\{
      \begin{array}{l}
        \left(1+\norma{\rho_2 (t)}_{\L1 (\reali^n;\reali^m)}\right)
        \norma{\rho_2 (t)}_{\L1 (\reali^n;\reali^m)} \,,
        \\
        \qquad\left(1+\norma{\rho_1 (t)}_{\L1 (\reali^n;\reali^m)}\right)
        \norma{\rho_1 (t)}_{\L1 (\reali^n;\reali^m)}
      \end{array}
    \right\}\d\tau
    \right)
    \\
    &
    & \;\times
      \left(\norma{\rho_o}_{\L\infty (\reali^n; \reali^m)} + \tv (\rho_o)\right) \;
      L_V \; \norma{\nabla \eta}_{\L\infty (\reali^n; \reali^n)} \;
      \norma{\rho_2 - \rho_1}_{\L1 ([0,t]\times\reali^n;\reali^m)}
    \\
    &
    & + \norma{\rho_o}_{\L\infty (\reali^n; \reali^m)} \;
      L_V \left(1 +
      \min \left\{\norma{\rho_1}_{\L1 (\reali^n; \reali^m)},
      \norma{\rho_2}_{\L1 (\reali^n; \reali^m)}\right\} \norma{\nabla \eta}_{\L\infty (\reali^n; \reali^n)}\right)
    \\
    &
    & \quad \times \norma{\Delta \eta}_{\L\infty (\reali^n; \reali)} \;
      \norma{\rho_2 - \rho_1}_{\L1 ([0,t]\times\reali^n;\reali^m)}
    \\
    & \leq
    & \exp
      \left(
      L_V \; \norma{D^2 \eta}_{\L\infty (\reali^n;\reali^{n\times n})} \;
      \norma{\rho_o}_{\L1 (\reali^n; \reali^m)} \; t
      \right)
    \\
    &
    & \quad\times
      \Bigl[
      1
      +
      L_V \left(\norma{D^2 \eta}_{\L\infty (\reali^n;\reali^{n\times n})}
      \norma{\Delta\eta}_{\L\infty (\reali^n; \reali)}
      + \norma{D^3 \eta}_{\L\infty (\reali^n;\reali^{n\times n\times n})}
      \right)
    \\
    &
    & \quad\qquad \times
      \left(1+\norma{\rho_o}_{\L1 (\reali^n;\reali^m)}\right)
      \norma{\rho_o}_{\L1 (\reali^n;\reali^m)} t
      \Bigr]
    \\
    &
    & \quad\times
      \left(\norma{\rho_o}_{\L\infty (\reali^n; \reali^m)} + \tv (\rho_o)\right) \;
      L_V \; \norma{\nabla \eta}_{\L\infty (\reali^n; \reali^n)} \;
      \norma{\rho_2 - \rho_1}_{\L1 ([0,t]\times\reali^n;\reali^m)}
    \\
    &
    & + \norma{\rho_o}_{\L\infty (\reali^n; \reali^m)} \;
      L_V
      \left(
      1+\norma{\rho_o}_{\L1 (\reali^n; \reali^m)}
      \norma{\nabla\eta}_{\L\infty (\reali^n; \reali^m)}
      \right)
      \norma{\Delta \eta}_{\L\infty (\reali^n; \reali)}    \\
    &
    & \qquad \times
      \norma{\rho_2 - \rho_1}_{\L1 ([0,t]\times\reali^n;\reali^m)}
    \\
    & \leq
    & C \;
      \left(\norma{\rho_o}_{\L\infty (\reali^n; \reali^m)} + \tv (\rho_o)\right)
      e^{C\, t} \,
      \norma{\rho_2 - \rho_1}_{\L1 ([0,t]\times\reali^n;\reali^m)}
  \end{eqnarray*}
  where $C$ is a positive constant depending on $\eta$, $L_V$ and
  $\norma{\rho_o}_{\L1 (\reali^n;\reali^m)}$. In particular, $C$ is
  independent of $\norma{\rho_o}_{\L\infty (\reali^n;\reali^m)}$,
  $\tv (\rho_o)$, $\rho_1$, $\rho_2$, $t$ and $T$.

  We thus obtain
  \begin{eqnarray*}
    &
    & \norma{\mathcal{T} (\rho_2) - \mathcal{T} (\rho_1)}_{\C0 ([0,T];\L1 (\reali^n; \reali^m))}
    \\
    & \leq
    & C \,
      \left(\norma{\rho_o}_{\L\infty (\reali^n; \reali^m)} + \tv (\rho_o)\right)
      e^{C\,T} \,
      \norma{\rho_2 - \rho_1}_{\C0 ([0,T];\L1 (\reali^n; \reali^m))} \,T
  \end{eqnarray*}
  proving that for
  $T \, e^{C\,T} < \frac1{C\left(\norma{\rho_o}_{\L\infty (\reali^n;
        \reali^m)} + \tv (\rho_o)\right)}$, $\mathcal{T}$ is a
  contraction and admits a unique fixed point $\rho^*$.  By
  Definition~\ref{def:solution}, $\rho^*$ solves on $[0,T]$ the Cauchy
  problem for~\eqref{eq:1} with initial datum $\rho_o$ assigned at
  time $t_o=0$. Moreover, $\rho^*$ satisfies
  $\norma{\rho^* (t)}_{\L1 (\reali^n; \reali^m)} = \norma{\rho_o}_{\L1
    (\reali^n; \reali^m)}$ for $t \in [0, T]$.

  \paragraph{Claim~3: Extension of $\rho^*$ for all times.} A repeated
  application of the previous step allows to iteratively construct
  maps
  $\rho^\nu \in \C0\left([T_{\nu-1},T_\nu]; \L1
    (\reali^n;\reali^m)\right)$ that solve
  \begin{displaymath}
    \left\{
      \begin{array}{l@{\qquad\quad}l}
        \partial_t \rho_i^\nu
        +
        \nabla \cdot
        \left( \rho^\nu_i \; V_i (\nabla \rho^\nu * \eta) \right)
        = 0
        & (t,x) \in  [T_{\nu-1},T_\nu] \times \reali^n
          \qquad i \in \{1, \ldots, m\}
        \\
        \rho^\nu (T_{\nu-1},x) = \rho^{\nu-1} (T_{\nu-1},x)
        & x \in \reali^n \,.
      \end{array}
    \right.
  \end{displaymath}
  The juxtaposition $\rho^*$ of $\rho_\nu$ for
  $\nu \in \naturali \setminus \{0\}$ solves~\eqref{eq:1} with initial
  datum $\rho_o$ assigned at time $t_o=0$ on the time interval
  $[0, \sup_\nu T_\nu [$.  At the same time, $\rho^*$ is a fixed point
  for $\mathcal{T}$ extended to the time interval
  $[0, \sup_\nu T_\nu [$.

  If $\sup_\nu T_\nu = +\infty$, then $\rho^*$ is globally defined on
  $\reali_+$. Otherwise, note that applying~\ref{item:3},
  \eqref{eq:18}, \eqref{eq:19}, we have for
  $t_1, t_2 \in \mathopen[0, \sup T_\nu\mathclose[$ with $t_1 < t_2$,
  \begin{eqnarray*}
    &
    & \norma{\rho^* (t_2) - \rho^* (t_1)}_{\L1 (\reali^n;\reali^m)}
    \\
    & =
    & \norma{
      \left(\Sigma (\Pi\rho^*)\right) (t_2)
      -
      \left(\Sigma (\Pi\rho^*)\right) (t_1)}_{\L1 (\reali^n;\reali^m)}
    \\
    & \leq
    & \left(
      1
      +
      \norma{\nabla (\nabla \cdot \Pi\rho^*)}_{\L1 ([0,t_2]\times \reali^n;\reali^{n\times m})}\right)
    \\
    &
    & \times
      \exp \norma{ \nabla \Pi\rho^*}_{\L1 ([0, t_2];\L\infty (\reali^n;\reali^{n^2\times m}))}
    \\
    &
    & \times
      \left(\norma{\rho_o}_{\L\infty (\reali^n; \reali^m)} + \tv (\rho_o)\right)
      \modulo{t_2 - t_1}
    \\
    & \leq
    & \left(
      1
      +
      L_V \left(
      \norma{D^2\eta}_{\L1 (\reali^n;\reali^{n\times n})}
      \norma{\Delta\eta}_{\L\infty (\reali^n;\reali)}
      +\norma{D^3\eta}_{\L1 (\reali^n;\reali^{n\times n\times n})}
      \right)
      \right.
    \\
    &
    & \quad\qquad\left. \times
      \left(1 + \norma{\rho_o}_{\L1 (\reali^n;\reali^m)}\right)
      \norma{\rho_o}_{\L1 (\reali^n;\reali^m)}
      t_2
      \right)
    \\
    &
    & \times
      \exp \left(
      L_V \, \norma{D^2\eta}_{\L\infty (\reali^n; \reali^{n\times n})} \,
      \norma{\rho_o}_{\L1 (\reali^n;\reali^m)}
      t_2
      \right)
    \\
    &
    & \times
      \left(\norma{\rho_o}_{\L\infty (\reali^n; \reali^m)} + \tv (\rho_o)\right)
      \modulo{t_2 - t_1}
    \\
    \\
    & \leq
    & \left(
      1
      +
      L_V \left(
      \norma{D^2\eta}_{\L1 (\reali^n;\reali^{n\times n})}
      \norma{\Delta\eta}_{\L\infty (\reali^n;\reali)}
      +\norma{D^3\eta}_{\L1 (\reali^n;\reali^{n\times n\times n})}
      \right)
      \right.
    \\
    &
    & \quad\qquad\left. \times
      \left(1 + \norma{\rho_o}_{\L1 (\reali^n;\reali^m)}\right)
      \norma{\rho_o}_{\L1 (\reali^n;\reali^m)}
      \sup_\nu T_\nu
      \right)
    \\
    &
    & \times
      \exp \left(
      L_V \, \norma{D^2\eta}_{\L\infty (\reali^n; \reali^{n\times n})} \,
      \norma{\rho_o}_{\L1 (\reali^n;\reali^m)}
      \sup_\nu T_\nu
      \right)
    \\
    &
    & \times
      \left(\norma{\rho_o}_{\L\infty (\reali^n; \reali^m)} + \tv (\rho_o)\right)
      \modulo{t_2 - t_1}
  \end{eqnarray*}
  which shows that $\rho^*$ is uniformly Lipschitz continuous on
  $[0, \sup_\nu T_\nu\mathclose[$. Thus, $\rho^* (\sup_\nu T_\nu)$ is
  well defined and can be used as initial datum for~\eqref{eq:1} to
  further extend $\rho^*$. This proves that $\rho^*$ can be extended
  to all $\reali_+$ and satisfies
  $\norma{\rho^* (t)}_{\L1 (\reali^n; \reali^m)} = \norma{\rho_o}_{\L1
    (\reali^n; \reali^m)}$ for $t \in \reali_+$.

  \paragraph{Claim~4: Uniqueness of the Solution to the Cauchy
    problem.} Let $\tilde\rho$ be a solution to~\eqref{eq:1} on the
  time interval $[0,\tilde T]$ in the sense of
  Definition~\ref{def:solution} with initial datum $\rho_o$ at time
  $t_o=0$. Then, for a $\bar\nu \in \naturali$,
  $\tilde T \in [T_{\bar \nu}, T_{\bar\nu+1}\mathclose[$, where we
  convene that $T_0 = 0$. By the uniqueness of the fixed point of a
  contraction, $\tilde \rho_{\strut\big| [T_\nu, T_{\nu+1}[}$
  coincides with $\rho^*_{\strut\big| [T_\nu, T_{\nu+1}[}$ for all
  $\nu = 0, \ldots, \bar\nu-1$. The same argument also applies to the
  time interval $[T_{\bar \nu}, \tilde T]$.

  \paragraph{Claim~5: Definition of the group $\mathcal{G}$.} By the
  arbitrariness of $\rho_o$, we can define a map
  \begin{displaymath}
    \mathcal{G}^+ \colon \reali_+ \times
    (\L1 \cap \L\infty \cap \BV) (\reali^n; \reali^m)
    \to
    (\L1 \cap \L\infty \cap \BV) (\reali^n; \reali^m)
  \end{displaymath}
  such that $t \mapsto \mathcal{G}^+_t \rho_o$ is the solution
  to~\eqref{eq:1} with datum $\rho_o$ assigned at time $0$ and
  satisfies
  $\norma{\mathcal{G}^+_t \rho_o}_{\L1 (\reali^n, \reali^m)} =
  \norma{\rho_o}_{\L1 (\reali^n; \reali^m)}$ for every
  $t \in \reali_+$.

  Note that the map $-V$ satisfies~\ref{item:7}, hence the above
  claims all apply to the equation
  \begin{equation}
    \label{eq:48}
    \partial_t \tilde\rho_i
    -
    \nabla \cdot \left( \tilde\rho_i \; V_i(\nabla \tilde\rho * \eta) \right)
    = 0
    \qquad i \in \{1, \ldots, m\}
  \end{equation}
  thus generating a map
  \begin{displaymath}
    \mathcal{G}^- \colon \reali_+ \times
    (\L1 \cap \L\infty \cap \BV) (\reali^n; \reali^m)
    \to
    (\L1 \cap \L\infty \cap \BV) (\reali^n; \reali^m)
  \end{displaymath}
  satisfying
  $\norma{\mathcal{G}^-_t \rho_o}_{\L1 (\reali^n, \reali^m)} =
  \norma{\rho_o}_{\L1 (\reali^n; \reali^m)}$ for all $t \in
  \reali_+$. Since~\eqref{eq:1} and~\eqref{eq:48} are autonomous and
  thanks to the representation~\eqref{eq:7}, the orbits of
  $\mathcal{G}^-$ are those of $\mathcal{G}^+$ but reversed in time,
  so that $\mathcal{G}^-_t (\mathcal{G}^+_t \rho_o) = \rho_o$, refer
  to Appendix~\ref{sec:time-reversibility} for more details. We thus
  define
  \begin{displaymath}
    \mathcal{G} (t,\rho_o) =
    \left\{
      \begin{array}{l@{\qquad}r@{\,}c@{\,}l}
        \mathcal{G}^- (-t,\rho_o)
        & t
        & <
        & 0
        \\
        \mathcal{G}^+ (t,\rho_o)
        & t
        & \geq
        & 0 \,.
      \end{array}
    \right.
  \end{displaymath}
  proving~\ref{item:33} and~\ref{item:32}. Note that $\mathcal{G}$
  satisfies
  $\norma{\mathcal{G} (t,\rho_o)}_{\L1 (\reali^n, \reali^m)} =
  \norma{\rho_o}_{\L1 (\reali^n; \reali^m)}$ for all $t \in \reali$,
  so that~\ref{item:constant_L1} is proved.

  \paragraph{Claim~6: Estimate~\ref{item:9}.}
  Fix
  $\widehat \rho_o, \widecheck \rho_o \in \left(\LL1 \cap \LL\infty
    \cap \BV\right) \left(\reali^n; \reali^m\right)$ and introduce the
  maps
  \begin{displaymath}
    \widehat{v} (t,x)
    =
    V\left(\left((\mathcal{G}_t \widehat \rho_o) * \eta\right) (x)\right)
    \quad\mbox{ and }\quad
    \widecheck{v} (t,x)
    =
    V\left(\left((\mathcal{G}_t \widecheck \rho_o) * \eta\right) (x)\right)
  \end{displaymath}
  as well as the solution $\bar\rho$ to the linear decoupled system
  \begin{displaymath}
    \left\{
      \begin{array}{l}
        \partial_t \bar\rho_i
        + \nabla \cdot
        \left( \bar\rho_i \;  \widecheck{v}_i (t,x)  \right)
        =0
        \qquad i\in\{1, \ldots, m\}
        \\
        \bar \rho (0) = \widehat \rho_o \,.
      \end{array}
    \right.
  \end{displaymath}
  Then,
  \begin{equation}
    \label{eq:12}
    \norma{\mathcal G_t \widehat{\rho}_o
      - \mathcal G_t\widecheck{\rho}_o }_{\L1 (\reali^n; \reali^m)}
    \leq
    \norma{\mathcal G_t\widehat{\rho}_o
      - \bar \rho (t)}_{\L1 (\reali^n; \reali^m)}
    +  \norma{\bar \rho (t) - \mathcal G_t\widecheck{\rho}_o}
    _{\L1 (\reali^n; \reali^m)}
  \end{equation}
  and we compute the two terms separately.  First, by~\ref{item:2},
  the first term in the right hand side of~\eqref{eq:12} can be
  bounded as follows:
  \begin{flalign*}
    & \norma{\mathcal G_t \widehat{\rho}_o - \bar \rho (t)}_{\L1
      (\reali^n; \reali^m)}
    \\
    \leq %
    & \exp \int_0^t \max\left\{ \norma{\nabla \widecheck v
        (\tau)}_{\L\infty (\reali^n; \reali^{n^2\times m})} ,
      \norma{\nabla \widehat v (\tau)}_{\L\infty (\reali^n;
        \reali^{n^2\times m})} \right\} \d\tau
    &[\mbox{By~\ref{item:2}}]
    \\
    & \quad\times \left( 1 {+} \int_0^t \max\left\{ \norma{\nabla
          \left(\nabla \cdot \widecheck v (\tau)\right)}_{\L1
          (\reali^n; \reali^{n\times m})} , \norma{\nabla \left(\nabla
            \cdot \widehat v (\tau)\right)}_{\L1 (\reali^n;
          \reali^{n\times m})} \right\}\d\tau \right)
    \\
    & \quad\times \left( \norma{\widehat \rho_o}_{\L\infty (\reali^n;
        \reali^m)} + \tv (\widehat \rho_o) \right) \; \norma{\widehat
      v - \widecheck v}_{\L1([0,t];\L\infty(\reali^n; \reali^{n\times
        m})}
    \\
    & + \norma{\widehat \rho_o}_{\L\infty (\reali^n; \reali^m)} \;
    \norma{\nabla \cdot \widehat v - \nabla\cdot \widecheck v}_{\L1
      ([0,t]\times \reali^n; \reali^m)}
    \\
    \leq %
    & \exp \left( L_V \norma{D^2 \eta}_{\L\infty (\reali^n;
        \reali^{n\times n})} \int_0^t \max \left\{ \norma{\mathcal
          G_\tau \widecheck{\rho}_o} _{\L1 (\reali^n;\reali^m)},
        \norma{\mathcal G_\tau \widehat{\rho}_o}_{\L1
          (\reali^n;\reali^m)} \right\}\d\tau \right) &
    [\mbox{By~\eqref{eq:18}}]
    \\
    & \quad\times \left( 1 + L_V \left(\norma{D^2\eta}_{\L1
          (\reali^n;\reali^{n\times n})} \norma{\Delta\eta}_{\L\infty
          (\reali^n; \reali)} + \norma{D^3\eta}_{\L1
          (\reali^n;\reali^{n\times n\times n})} \right) \right.
    \\
    & \qquad \times \int_0^t \max\left\{ \left(1+\norma{\mathcal
          G_\tau \widecheck{\rho}_o}_{\L1 (\reali^n;\reali^m)}\right)
      \norma{\mathcal G_\tau\widecheck{\rho}_o} _{\L1
        (\reali^n;\reali^m)}, \right.
    \\
    &\qquad\qquad\qquad \left.  \left.  \left(1+\norma{\mathcal G_\tau
            \widehat{\rho}_o }_{\L1 (\reali^n;\reali^m)}\right)
        \norma{\mathcal G_\tau \widehat{\rho}_o}_{\L1
          (\reali^n;\reali^m)} \right\}\d\tau \right) &
    [\mbox{By~\eqref{eq:19}}]
    \\
    & \quad\times \left( \norma{\widehat \rho_o}_{\L\infty (\reali^n;
        \reali^m)} {+} \tv (\widehat \rho_o) \right) L_V
    \norma{\nabla\eta}_{\L\infty (\reali^n;\reali^n)} \norma{\mathcal
      G \widecheck{\rho}_o {-} \mathcal G
      \widehat{\rho}_o}_{\L1([0,t];\L1(\reali^n; \reali^m))} &
    [\mbox{By~\eqref{eq:20}}]
    \\
    & + \norma{\widehat \rho_o}_{\L\infty (\reali^n; \reali^m)}
    \\
    & \quad \times L_V \left(1 + \min \left\{\norma{\mathcal G
          \widehat{\rho}_o}_{\L1 ([0,t]\times \reali^n; \reali^m)},
        \norma{\mathcal G\widecheck{\rho}_o}_{\L1
          ([0,t]\times\reali^n; \reali^m)}\right\} \norma{\nabla
        \eta}_{\L1 (\reali^n; \reali^n)}\right)
    \\
    & \quad \times \norma{\Delta \eta}_{\L1 (\reali^n; \reali)}
    \norma{\mathcal G\widehat{\rho}_o - \mathcal
      G\widecheck{\rho}_o}_{\L1 ([0,t];\L1 (\reali^n; \reali^m))} &
    [\mbox{By~\eqref{eq:22}}]
    \\
    \leq %
    & \exp \left( L_V \norma{D^2 \eta}_{\L\infty (\reali^n;
        \reali^{n\times n})} \; \max \left\{
        \norma{\widecheck\rho_o}_{\L1 (\reali^n;\reali^m)},
        \norma{\widehat\rho_o}_{\L1 (\reali^n;\reali^m)} \right\}\, t
    \right) & [\mbox{By~\eqref{eq:13}}]
    \\
    & \quad\times \left( 1 + L_V \left(\norma{D^2\eta}_{\L1
          (\reali^n;\reali^{n\times n})} \norma{\Delta\eta}_{\L\infty
          (\reali^n; \reali)} + \norma{D^3\eta}_{\L1
          (\reali^n;\reali^{n\times n\times n})} \right) \right.
    \\
    & \qquad\qquad \times \max\left\{
      \left(1+\norma{\widecheck\rho_o}_{\L1
          (\reali^n;\reali^m)}\right) \norma{\widecheck\rho_o}_{\L1
        (\reali^n;\reali^m)}, \right.
    \\
    & \qquad\qquad\qquad\qquad\qquad
    \left. \left. \left(1+\norma{\widehat\rho_o}_{\L1
            (\reali^n;\reali^m)}\right) \norma{\widehat\rho_o}_{\L1
          (\reali^n;\reali^m)} \right\} t \right) &
    [\mbox{By~\eqref{eq:13}}]
    \\
    & \quad\times \left( \norma{\widehat \rho_o}_{\L\infty (\reali^n;
        \reali^m)} {+} \tv (\widehat \rho_o) \right) L_V
    \norma{\nabla\eta}_{\L\infty (\reali^n;\reali^n)} \norma{\mathcal
      G\widecheck{\rho}_o {-} \mathcal G\widehat{\rho}_o
    }_{\L1([0,t];\L1(\reali^n; \reali^m))}
    \\
    & + \norma{\widehat \rho_o}_{\L\infty (\reali^n; \reali^m)} L_V
    \left(1 {+} \min \left\{\norma{\widehat \rho_o}_{\L1 (\reali^n;
          \reali^m)}, \norma{\widecheck \rho_o}_{\L1 (\reali^n;
          \reali^m)}\right\} \, t \, \norma{\nabla \eta}_{\L1
        (\reali^n; \reali^n)}\right) & [\mbox{By~\eqref{eq:13}}]
    \\
    & \quad \times \norma{\Delta \eta}_{\L1 (\reali^n; \reali)}
    \norma{\mathcal G\widehat{\rho}_o - \mathcal
      G\widecheck{\rho}_o}_{\L1 ([0,t];\L1 (\reali^n; \reali^m))} \,.
  \end{flalign*}
  Introducing a positive, continuous and non decreasing map
  $C = C (t)$ dependent on $\eta$ and on the initial data
  $\widecheck \rho_o$, $\widehat \rho_o$, we obtain a bound of the
  form, for all $t \in \reali_+$,
  \begin{equation}
    \label{eq:24}
    \norma{\mathcal G_t \widehat{\rho}_o - \bar \rho (t)}_{\L1 (\reali^n; \reali^m)}
    \leq
    C(t) \,
    \norma{\mathcal G \widehat{\rho}_o -
      \mathcal G \widecheck{\rho}_o}_{\L1 ([0,t];\L1 (\reali^n; \reali^m))} \,.
  \end{equation}

  The second term in the right hand side in~\eqref{eq:12} is estimated
  by means of~\ref{item:13}:
  \begin{equation}
    \label{eq:21}
    \norma{\bar \rho (t) - \mathcal G_t \widecheck{\rho}_o}_{\L1 (\reali^n; \reali^m)}
    \leq
    \norma{\widecheck \rho_o - \widehat \rho_o}_{\L1 (\reali^n; \reali^m)} \,.
  \end{equation}
  Grouping~\eqref{eq:24} and~\eqref{eq:21} we get
  \begin{displaymath}
    \norma{\mathcal G_t \widecheck{\rho}_o
      - \mathcal G_t \widehat{\rho}_o}_{\L1 (\reali^n; \reali^m)}
    \leq
    \norma{\widecheck \rho_o - \widehat \rho_o}_{\L1 (\reali^n; \reali^m)}
    + C (t) \, \norma{\mathcal G \widehat{\rho}_o -
      \mathcal G \widecheck{\rho}_o}_{\L1 ([0,t];\L1 (\reali^n; \reali^m))}
  \end{displaymath}
  so that, after renaming the constant $C$, an application of Gronwall
  Lemma yields the desired estimate in~\ref{item:9}.

  \paragraph{Claim 7: Stability Estimate~\ref{item:14}.}
  Assume for simplicity $t \geq 0$, the other case being
  analogous. Set
  \begin{displaymath}
    \widehat v (t,x)
    =
    \widehat V\left(\nabla
      (\widehat {\mathcal G}_t \rho_o * \eta) (x)
    \right)
    \quad \mbox{ and } \quad
    \widecheck v (t,x)
    =
    \widecheck V\left(\nabla
      (\widecheck {\mathcal G}_t\rho_o * \eta) (x)
    \right) \,.
  \end{displaymath}
  In view of the application of~\ref{item:2}, call
  $B = B_{\reali^{n\times m}} (0, \norma{\rho_o}_{\L1 (\reali^n;
    \reali^m)} \, \norma{\nabla\eta}_{\L\infty (\reali^n;
    \reali^n)})$, note that
  $\norma{\nabla\left(\widecheck{\mathcal G}_t \rho_o * \eta
    \right)}_{\LL\infty\left(\reali^n; \reali^{n \times m}\right)} \le
  \norma{\rho_o}_{\LL1\left(\reali^n; \reali^m\right)} \norma{\nabla
    \eta}_{\LL\infty\left(\reali^n; \reali^n\right)}$ for every
  $t \in \reali$. Compute for $t \in \reali$,
  \begin{eqnarray*}
    &
    & \norma{\widehat{v} (t) - \widecheck{v} (t)}_{\L\infty (\reali^n; \reali^{n\times m})}
    \\
    & \leq
    & \norma{\widehat{v} (t) -
      \widehat V\left(\nabla (\widecheck {\mathcal G}_t
      \rho_o * \eta)\right)}_{\L\infty (\reali^n; \reali^{n\times m})}
      +
      \norma{\widehat V\left(\nabla (\widecheck {\mathcal G}_t
      \rho_o * \eta)\right) - \widecheck{v} (t)}_{\L\infty (\reali^n; \reali^{n\times m})}
    \\
    & \leq
    & L_V \norma{
      \nabla
      (\widehat {\mathcal G}_t \rho_o * \eta)
      -
      \nabla (\widecheck {\mathcal G}_t \rho_o * \eta)
      }_{\L\infty (\reali^n;\reali^{n\times m})}
      +
      \norma{\widehat{V} - \widecheck{V}}_{\L\infty (B;\reali^{n\times m})}
    \\
    & \leq
    & L_V \norma{
      \widehat {\mathcal G}_t \rho_o
      -
      \widecheck {\mathcal G}_t \rho_o
      }_{\L1 (\reali^n;\reali^{n\times m})}
      \norma{\nabla\eta}_{\L\infty (\reali^n; \reali^n)}
      +
      \norma{\widehat{V} - \widecheck{V}}_{\L\infty (B;\reali^{n\times m})}
  \end{eqnarray*}
  so that
  \begin{eqnarray}
    \nonumber
    &
    & \norma{\widehat{v} - \widecheck{v}}_{\L1([0,t];\L\infty (\reali^n; \reali^{n\times m}))}
    \\
    \label{eq:999}
    & \leq
    & L_V \norma{
      \widehat {\mathcal G} \rho_o
      -
      \widecheck {\mathcal G} \rho_o
      }_{\L1 ([0,t]\times\reali^n;\reali^{n\times m})}
      \norma{\nabla\eta}_{\L\infty (\reali^n; \reali^n)}
      +
      \norma{\widehat{V} - \widecheck{V}}_{\L\infty (B;\reali^{n\times m})} t \,.
  \end{eqnarray}
  Similarly,
  \begin{eqnarray*}
    &
    & \norma{
      \nabla \cdot \widehat{v} (t) - \nabla \cdot \widecheck{v} (t)
      }_{\L1 (\reali^n; \reali^m)}
    \\
    & =
    & \norma{
      \nabla \cdot \left(\widehat{v} (t) - \widecheck{v} (t)\right)
      }_{\L1 (\reali^n; \reali^m)}
    \\
    & \leq
    & \norma{D
      \left(\widehat V\left(\nabla
      (\widehat {\mathcal G}_t \rho_o * \eta)
      \right)
      -
      \widecheck V\left(\nabla
      (\widecheck {\mathcal G}_t\rho_o * \eta)
      \right)\right)}_{\L1 (\reali^n; \reali^{(n\times m)^2})}
    \\
    &
    & \quad \times \norma{
      D^2 \left((\widehat {\mathcal G}_t \rho_o * \eta)
      -
      (\widecheck {\mathcal G}_t\rho_o * \eta)
      \right)}_{\L\infty (\reali^n; \reali^{n^2 \times m})}
    \\
    & \leq
    & \Bigg[
      \norma{D
      \left(\widehat V\left(\nabla
      (\widehat {\mathcal G}_t \rho_o * \eta)
      \right)
      -
      \widehat V\left(\nabla
      (\widecheck {\mathcal G}_t\rho_o * \eta)
      \right)\right)}_{\L1 (\reali^n; \reali^{(n\times m)^2})}
    \\
    &
    & \qquad +
      \norma{D
      \left(\widehat V\left(\nabla
      (\widecheck {\mathcal G}_t \rho_o * \eta)
      \right)
      -
      \widecheck V\left(\nabla
      (\widecheck {\mathcal G}_t\rho_o * \eta)
      \right)\right)}_{\L1 (\reali^n; \reali^{(n\times m)^2})}
      \Bigg]
    \\
    &
    & \quad \times \norma{
      D^2 \left((\widehat {\mathcal G}_t \rho_o -
      \widecheck {\mathcal G}_t\rho_o )* \eta
      \right)}_{\L\infty (\reali^n; \reali^{n^2 \times m})}
    \\
    & \leq
    & \left[
      L_V \norma{
      \nabla
      (\widehat {\mathcal G}_t \rho_o * \eta)
      -
      \nabla (\widecheck {\mathcal G}_t \rho_o * \eta)
      }_{\L1 (\reali^n;\reali^{n\times m})}
      +
      \norma{D\widehat{V} - D\widecheck{V}}_{\L1
      (B;\reali^{\left(n\times m\right)^2})}
      \right]
    \\
    &
    & \quad \times
      \norma{
      \widehat {\mathcal G}_t \rho_o
      -
      \widecheck {\mathcal G}_t \rho_o
      }_{\L1 (\reali^n;\reali^{n\times m})}
      \norma{D^2\eta}_{\L\infty (\reali^n; \reali^{n\times n})}
    \\
    & \leq
    &       \left[
      L_V \norma{
      \widehat {\mathcal G}_t \rho_o
      -
      \widecheck {\mathcal G}_t \rho_o
      }_{\L1 (\reali^n;\reali^{n\times m})}
      \norma{\nabla\eta}_{\L1 (\reali^n; \reali^n)}
      +
      \norma{D\widehat{V} - D\widecheck{V}}_{\L1 (B;\reali^{\left(n\times m\right)^2})}
      \right]
    \\
    &
    & \quad \times
      \norma{
      \widehat {\mathcal G}_t \rho_o
      -
      \widecheck {\mathcal G}_t \rho_o
      }_{\L1 (\reali^n;\reali^{n\times m})}
      \norma{D^2\eta}_{\L\infty (\reali^n; \reali^{n\times n})}
  \end{eqnarray*}
  so that
  \begin{eqnarray}
    \nonumber
    &
    & \norma{
      \nabla \cdot \widehat{v} - \nabla \cdot \widecheck{v}
      }_{\L1 ([0, t] \times \reali^n; \reali^m)}
    \\
    \nonumber
    & \le
    & \left[
      L_V \norma{
      \widehat {\mathcal G} \rho_o
      -
      \widecheck {\mathcal G} \rho_o
      }_{\L1 ([0, t] \times \reali^n;\reali^{n\times m})}
      \norma{\nabla\eta}_{\L1 (\reali^n; \reali^n)}
      +
      \norma{D\widehat{V} - D\widecheck{V}}_{\L1 (B;\reali^{\left(n\times m\right)^2})} t
      \right]
    \\
    &
    & \quad \times 2 \norma{\rho_o}_{\LL1\left(\reali^n; \reali^m\right)}
      \norma{D^2\eta}_{\L\infty (\reali^n; \reali^{n\times n})}.
      \label{eq:998}
  \end{eqnarray}
  Therefore we deduce, for every $t \in \reali$,
  \begin{eqnarray*}
    &
    & \norma{ \widehat {\mathcal G}_t \rho_o - \widecheck {\mathcal
      G}_t \rho_o}_{\L1 (\reali^n; \reali^m)}
    \\
    & \leq
    & \exp \int_0^t \max\left\{ \norma{\nabla \widecheck v
      (\tau)}_{\L\infty (\reali^n; \reali^{n^2\times m})} ,
      \norma{\nabla \widehat v (\tau)}_{\L\infty (\reali^n;
      \reali^{n^2\times m})} \right\} \d\tau %
      \qquad\qquad\qquad\;  [\mbox{By~\ref{item:2}}]
    \\
    &
    & \quad\times \left( 1 + \int_0^t \max\left\{ \norma{\nabla
      \left(\nabla \cdot \widecheck v (\tau)\right)}_{\L1
      (\reali^n; \reali^{n\times m})} , \norma{\nabla \left(\nabla
      \cdot \widehat v (\tau)\right)}_{\L1 (\reali^n;
      \reali^{n\times m})} \right\}\d\tau \right)
    \\
    &
    & \quad\times \left( \norma{\rho_o}_{\L\infty (\reali^n;
      \reali^m)} + \tv (\rho_o) \right) \norma{\widehat v -
      \widecheck v}_{\L1([0,t];\L\infty(\reali^n; \reali^{n\times
      m}))}
    \\
    &
    & + \norma{\rho_o}_{\L\infty (\reali^n; \reali^m)} \;
      \norma{\nabla \cdot \widehat v - \nabla\cdot \widecheck v}_{\L1
      ([0,t]\times \reali^n; \reali^m)}
    \\
    & \leq %
    & \exp \int_0^t \max\left\{ \norma{\nabla \widecheck v
      (\tau)}_{\L\infty (\reali^n; \reali^{n^2\times m})} ,
      \norma{\nabla \widehat v (\tau)}_{\L\infty (\reali^n;
      \reali^{n^2\times m})} \right\} \d\tau %
      \qquad\quad [\mbox{By~\eqref{eq:999}, \eqref{eq:998}}]
    \\
    &
    & \quad\times \left( 1 + \int_0^t \max\left\{ \norma{\nabla
      \left(\nabla \cdot \widecheck v (\tau)\right)}_{\L1
      (\reali^n; \reali^{n\times m})} , \norma{\nabla \left(\nabla
      \cdot \widehat v (\tau)\right)}_{\L1 (\reali^n;
      \reali^{n\times m})} \right\}\d\tau \right)
    \\
    &
    & \quad\times \left( \norma{\rho_o}_{\L\infty (\reali^n;
      \reali^m)} + \tv (\rho_o) \right)
    \\
    &
    & \quad \times \left[ L_V \norma{ \widehat {\mathcal G} \rho_o -
      \widecheck {\mathcal G} \rho_o }_{\L1
      ([0,t]\times\reali^n;\reali^{n\times m})}
      \norma{\nabla\eta}_{\L1 (\reali^n; \reali^n)} +
      \norma{\widehat{V} - \widecheck{V}}_{\L\infty (B;\reali^{n\times
      m})} t \right]
    \\
    &
    & + 2 \norma{\rho_o}^2_{\L\infty (\reali^n; \reali^m)} \norma{D^2
      \eta}_{ \LL\infty\left(\reali^n; \reali^{n \times n}\right)}
    \\
    &
    & \quad \times \left[ L_V \norma{ \widehat {\mathcal G} \rho_o {-}
      \widecheck {\mathcal G} \rho_o }_{\L1 ([0, t] \times
      \reali^n;\reali^{n\times m})} \norma{\nabla\eta}_{\L1
      (\reali^n; \reali^n)} {+} \norma{D\widehat{V} {-}
      D\widecheck{V}}_{\L1 (B;\reali^{\left(n\times m\right)^2})} t
      \right] \,
  \end{eqnarray*}
  so that
  \begin{eqnarray*}
    \norma{ \widehat {\mathcal G}_t \rho_o - \widecheck {\mathcal
    G}_t \rho_o}_{\L1 (\reali^n; \reali^m)}
    & \leq
    & C \, e^{C\, t} \,t
      \left(
      \norma{\widehat{V} - \widecheck{V}}_{\L\infty (B;\reali^{n\times m})}
      +
      \norma{D\widehat{V} - D\widecheck{V}}_{\L1 (B;\reali^{\left(n\times m\right)^2})}
      \right)
    \\
    &
    & + C\, \norma{
      \widehat {\mathcal G} \rho_o
      -
      \widecheck {\mathcal G} \rho_o
      }_{\L1 ([0, t] \times \reali^n;\reali^{n\times m})}
  \end{eqnarray*}
  An application of Gronwall Lemma completes the proof of the claim.

  \paragraph{Regularity and Positivity:} The statements~\ref{item:12}
  and~\ref{item:4} follow from formula~\eqref{eq:7} with $v_i$ as
  in~\eqref{eq:25}.

  The proof of Theorem~\ref{thm:cauchyProblem} is completed.
\end{proofof}

\begin{proofof}{Proposition~\ref{prop:polar}}
  Denote $\rho (t,x) = (\mathcal{G}_t\rho_o) (x)$ and define
  $r_o (x) = \rho_o (R\, x)$, $r (t,x) = (\mathcal{G}_t r_o) (x)$.

  For later use, observe that if $\psi \in \Cc1 (\reali^2;\reali)$,
  setting $\phi (t,x) = \psi (t, R^{-1}x)$, we then have
  \begin{equation}
    \label{eq:27}
    \partial_t \phi (t, x) = \partial_t \psi (t, R^{-1}x)
    \mbox{ and }
    \nabla\phi (t,x) \; R= \nabla\psi (t, R^{-1}x)\,.
  \end{equation}
  Moreover, for $i\in \{1, \ldots, m\}$,
  \begin{flalign}
    \nonumber \left(\nabla r_i (t) *\eta\right) (x) = & \nabla
    \int_{\reali^n} r_i (t,\xi) \, \eta (x-\xi) \d\xi
    \\
    \nonumber = & \int_{\reali^n} r_i (t,\xi) \, \nabla \eta (x-\xi)
    \d\xi
    \\
    \nonumber = & \int_{\reali^n} \rho_i (t,R\, \xi) \, \nabla \eta
    (x-\xi) \d\xi
    \\
    \nonumber = & \int_{\reali^n} \rho_i (t,z) \, \nabla \eta
    (x-R^{-1} z) \d{z}
    \\
    \nonumber = & \int_{\reali^n} \rho_i (t,z) \, \nabla \eta (R\,x-
    z) \d{z} \, R &[\mbox{Since } \nabla\eta (x) = \nabla\eta
    (R\,x)\,R \mbox{ by~\ref{item:18}}]
    \\
    \label{eq:49}
    = & \left(\nabla \rho_i (t) *\eta\right) (R\, x) \, R \,.
  \end{flalign}
  Using~\ref{item:19}, \ref{item:18}, \ref{item:20}, \eqref{eq:27}
  and~\eqref{eq:49}, we now verify that $r$ solves~\eqref{eq:1} with
  initial datum $r_o$.
  \begin{eqnarray*}
    &
    & \int_{\reali_+} \int_{\reali}
      r_i (t,x) \, \partial_t \psi (t,x)
      \d{x} \d{t}
    \\
    &
    & + \int_{\reali_+} \int_{\reali}
      r_i (t,x)\, \nabla \psi (t,x) \,
      V_i\left(\left(\nabla r_1 (t) *\eta\right) (x)
      , \ldots,
      \left(\nabla r_m (t) *\eta\right) (x)\right)
      \d{x} \d{t}
    \\
    &
    & - \int_{\reali} r_{o,i} (t,x) \, \psi (0,x) \d{x}
    \\
    & =
    & \int_{\reali_+} \int_{\reali}
      \rho^*_i (t,R\,x) \, \partial_t \psi (t,x)
      \d{x} \d{t}
    \\
    &
    & + \int_{\reali_+} \int_{\reali}
      \rho^*_i (t,R\,x)\, \nabla \psi (t,x) \,
      V_i\left(\left(\nabla \rho_1 (t) *\eta\right) (R\, x) \, R
      , \ldots,
      \left(\nabla \rho_m (t) *\eta\right) (R\, x) \, R\right)
      \d{x} \d{t}
    \\
    &
    & - \int_{\reali} \rho_{o,i}^* (t,R\,x) \, \psi (0,x) \d{x}
    \\
    & =
    & \int_{\reali_+} \int_{\reali}
      \rho^*_i (t,y) \, \partial_t \psi (t,R^{-1}y)
      \d{y} \d{t}
    \\
    &
    & + \int_{\reali_+} \int_{\reali}
      \rho^*_i (t,y)\, \nabla \psi (t,R^{-1}y) \,
      V_i\left(\left(\nabla \rho_1 (t) *\eta\right) (y) \, R
      , \ldots,
      \left(\nabla \rho_m (t) *\eta\right) (y) \, R\right)
      \d{y} \d{t}
    \\
    &
    & - \int_{\reali} \rho_{o,i}^* (t,y) \, \psi (0,R^{-1}y) \d{y}
    \\
    & =
    & \int_{\reali_+} \int_{\reali}
      \rho^*_i (t,y) \, \partial_t \phi (t,y)
      \d{y} \d{t}
    \\
    &
    & + \int_{\reali_+} \int_{\reali}
      \rho^*_i (t,y)\, \nabla \phi (t,y) \, R
      \, R^{-1} \, V_i\left(\left(\nabla \rho_1 (t) *\eta\right) (y)
      , \ldots,
      \left(\nabla \rho_m (t) *\eta\right) (y) \right)
      \d{y} \d{t}
    \\
    &
    & - \int_{\reali} \rho_{o,i}^* (t,y) \, \phi (0,y) \d{y}
      =
      0
  \end{eqnarray*}
  The uniqueness proved in Theorem~\ref{thm:cauchyProblem} thus
  implies that $\rho^* = r$, completing the proof.
\end{proofof}

Observe for later use that the formula~\eqref{eq:7} implies that, for
every $i \in \left\{1, \ldots, m\right\}$,
\begin{equation}
  \label{eq:39}
  \spt \rho_i (t)
  = X_i
  \left(t; 0, \spt \rho_{o,i} \right) \,.
\end{equation}

\begin{proofof}{Lemma~\ref{lem:easy}}
  Using~\eqref{eq:29} and~\eqref{eq:25}, we obtain, for
  $i\in \{1, \ldots, m\}$, the estimate
  \begin{flalign*}
    \norma{\dot X_i (t; 0, x_o)} = %
    & \norma{V_i\left(\nabla\left(\rho (t) * \eta\right) (x)\right)} %
    & [\mbox{By~\eqref{eq:25}--\eqref{eq:29}}]
    \\
    \leq & %
    L_V \, \norma{\nabla\left(\rho (t) * \eta\right)}
    _{\L\infty(\reali^n; \reali^{n\times m})} %
    & [\mbox{By~\ref{item:7}}]
    \\
    \leq & %
    L_V \; \norma{\nabla\eta}_{\L\infty (\reali^n; \reali^n)} \;
    \norma{\rho_o}_{\L1 (\reali^n; \reali^m)}\,, %
    & [\mbox{By~\eqref{eq:16} and~\ref{item:constant_L1} in
      Theorem~\ref{thm:cauchyProblem}}]
  \end{flalign*}
  the latter quantity being $W$ as in~\eqref{eq:9}. Hence,
  by~\eqref{eq:39}, \eqref{eq:7} and the assumption on $\spt\rho_o$,
  \begin{displaymath}
    \spt \rho (t)
    \subseteq
    \bigcup_{i=1}^m \spt X_i\left(t; 0, \spt \rho_{o,i}\right)
    \subseteq
    \bigcup_{i=1}^m \spt X_i\left(t; 0, B (x_o, r)\right)
    \subseteq B\left(x_o, r + W\,t\right) \,,
  \end{displaymath}
  completing the proof.
\end{proofof}

\begin{proofof}{Proposition~\ref{prop:Caratt}}
  Consider first~\ref{item:22}.  By contradiction, referring
  to~\eqref{eq:39}, assume there exist
  $\bar t \in \mathopen]0, +\infty\mathclose[$ and
  $x_o \in \spt \rho_{o,i}$ such that
  $X_i (\bar t; 0, x_o) \in \partial C$,
  $X_i \left([0, \bar t\mathclose[;0, \spt \rho_{o,i} \right)
  \subseteq C$ and for a suitable $\epsilon > 0$,
  $X_i (\mathopen]\bar t, \bar t + \epsilon\mathclose[, 0, x_o)
  \subset \reali^n \setminus C$. Call $\bar x = X_i (\bar t; 0,
  x_o)$. Since $C$ is non empty and closed, there exists a vector
  $\nu$ in the (outer) normal cone~\cite[\S~1.2,
  Formula~(6)]{MR709590} to $C$ at $\bar x$ such that
  \begin{equation}
    \label{eq:31}
    V_i\left(\nabla \left(\rho (\bar t) * \eta\right) (\bar x) \right)
    \cdot \nu \geq 0 \,.
  \end{equation}
  On the other hand, using~\ref{item:17}, \ref{item:21}, \ref{item:24}
  and the convexity of $C$,
  \begin{eqnarray}
    \label{eq:33}
    V_i\left(\nabla \left(\rho (\bar t) * \eta\right) (\bar x) \right)
    \cdot \nu
    & =
    & v\left(\norma{\nabla \left(\rho (\bar t) * \eta\right) (\bar x)}\right) \;
      \nabla \left(\rho_i (\bar t) * \eta\right) (\bar x)
      \cdot \nu \,;
    \\
    \label{eq:34}
    \left(\nabla \left(\rho_i (\bar t) * \eta\right) (\bar x)\right) \cdot
    \nu
    & =
    & \left(\left(\rho_i (\bar t) * \nabla \eta\right) (\bar x)\right) \cdot \nu
    \\
    \nonumber
    & =
    & \int_{\reali^n} \rho_i (\bar t,\xi) \; \nabla\eta (\bar x-\xi) \cdot \nu \d\xi
    \\
    \nonumber
    & =
    & \int_{\spt \rho_i (\bar t)\cap B (\bar x,\ell)}
      \rho_i (\bar t,\xi) \; \nabla\eta (\bar x-\xi) \cdot \nu \d\xi
    \\
    \label{eq:35}
    & =
    & \int_{\spt \rho_i (\bar t)\cap B (\bar x,\ell)}
      \underbrace{\rho_i (\bar t,\xi)}_{> 0\mbox{ a.e.}} \;
      \underbrace{\tilde\eta' \left(\norma{\bar x-\xi}\right)}_{< 0\mbox{ a.e.}} \;
      \underbrace{\frac{\bar x - \xi}{\norma{\bar x - \xi}} \cdot \nu}_{> 0\mbox{ a.e.}} \d\xi
  \end{eqnarray}
  hence
  \begin{displaymath}
    V_i\left(\nabla \left(q (\bar t) * \eta\right) (\bar x) \right)
    \cdot \nu \quad
    \left\{
      \begin{array}{l@{\quad\mbox{ if }}l}
        <0
        & \mathcal{L}\left(\spt \rho_i (\bar t) \cap B (\bar x,\ell)\right) >0
        \\
        =0
        &  \mathcal{L}\left(\spt \rho_i (\bar t) \cap B (\bar x,\ell)\right) =0
      \end{array}
    \right.
  \end{displaymath}
  contradicting~\eqref{eq:31} in case
  $\mathcal{L}\left(\spt \rho_i (\bar t) \cap B (\bar x,\ell)\right)
  >0$. On the other hand, by~\eqref{eq:41},
  \begin{displaymath}
    \mathcal{L}\left(\spt \rho_i (\bar t) \cap B (\bar x,\ell)\right) =
    0
    \implies
    \rho_i (\bar t) = 0 \mbox{ a.e.~on } B (\bar x, \ell)
    \implies
    B (\bar x, \ell) \subseteq \left(\reali^n \setminus \spt \rho_i (t)\right)
  \end{displaymath}
  where we used the fact that $B (\bar x,\ell)$ is an open set. The
  latter statement contradicts the choice of $\bar x$. Indeed,
  by~\eqref{eq:39},
  \begin{displaymath}
    \bar x
    = X_i (\bar t; 0, x_o)
    \in X_i \left(\bar t; 0, \spt \rho_{o,i} \right)
    = \spt \rho_i (t) \,,
  \end{displaymath}
  completing the proof of~\ref{item:22}. Now, \ref{item:23}
  immediately follows.
\end{proofof}

\begin{proofof}{Corollary~\ref{cor:clustering}}
  First, we prove that if $x_o \in B (x_h, r_h)$, the curve
  $t \mapsto X_i (t;0,x_o)$ defined in~\eqref{eq:29} remains in the
  same ball $B (x_h, r_h)$ for all times $t\geq 0$. By contradiction,
  assume there exists a first time
  $\bar t \in \mathopen]0, +\infty\mathclose[$ when the characteristic
  $t \mapsto X_i (t;0,x_o)$ exits the ball $B (x_h, r_h)$. Then,
  $X_i (\bar t; 0, x_o) \in \partial B (x_h, r_h)$,
  $X_i \left([0, \bar t\mathclose[;0, x_o \right) \subseteq B (x_h,
  r_h)$ and for a suitable $\epsilon > 0$,
  $X_i (\mathopen]\bar t, \bar t + \epsilon\mathclose[; 0, x_o) \not
  \in B (x_h, r_h)$.

  Call $\bar x = X_i (\bar t; 0, x_o)$. Fix an exterior normal $\nu$
  to $B (x_h, r_h)$ at $\bar x$, use equality~\eqref{eq:33} and repeat
  computations analogous to those between~\eqref{eq:34}
  and~\eqref{eq:35} to obtain:
  \begin{displaymath}
    \left(\nabla \left(\rho_i (\bar t) * \eta\right) (\bar x)\right)
    \cdot \nu
    =
    \int_{B (x_h, r_h) \cap B (\bar x,\ell)}
    \underbrace{\rho_i (\bar t,\xi)}_{> 0\mbox{ a.e.}} \;
    \underbrace{\tilde\eta' \left(\norma{\bar x-\xi}\right)}_{< 0\mbox{ a.e.}} \;
    \underbrace{\frac{\bar x - \xi}{\norma{\bar x-\xi}} \cdot \nu}_{> 0\mbox{ a.e.}} \, \d\xi
    <
    0
  \end{displaymath}
  since by construction
  $\mathcal{L}\left(B (x_{\bar h}, r_{\bar h}) \cap B (\bar
    x,\ell)\right) > 0$, getting a contradiction.

  From expression~\eqref{eq:7} the inclusion~\eqref{eq:28} then
  follows. Equality~\eqref{eq:26} is a consequence of the conservative
  form of~\eqref{eq:1}.
\end{proofof}

\begin{proofof}{Proposition~\ref{prop:stationary}}
  We verify that the map $\rho$ given by $\rho (t,x) = \rho_o (x)$
  solves~\eqref{eq:1} in the sense of
  Definition~\ref{def:solution}. To this aim, fix $x \in \spt\rho_o$
  so that $x \in B (x_{\bar h},r/2)$ for a suitable
  $\bar h \in \{1, \ldots,k\}$. Note that for
  $i \in \{1, \ldots, m\}$,
  \begin{flalign*}
    \left(\nabla \rho_i (t) * \eta\right) (x) =%
    & \int_{\reali^n} (\rho_o)_i (y) \; \nabla\eta (x-y) \d{y}%
    & [\mbox{Since } \rho (t) = \rho_o]
    \\
    = & \int_{\bigcup_{h=1}^k B (x_h,r / 2)} (\rho_o)_i (y) \;
    \nabla\eta (x-y) \d{y}%
    &[\mbox{By~\ref{item:31}}]
    \\
    = & \sum_{h=1}^k \int_{B (x_h,r / 2)} (\rho_o)_i (y) \; \nabla\eta
    (x-y) \d{y} \,.%
    &[\mbox{By~\ref{item:31}}]
  \end{flalign*}
  Each term in the latter sum vanishes by~\ref{item:30}. Indeed, when
  $h=\bar h$
  \begin{displaymath}
    \begin{array}{rcl}
      x \in B (x_{\bar h}, r/2)
      & \implies
      & \norma{x - x_{\bar h}} \leq r/2
      \\
      y \in B (x_{\bar h}, r/2)
      & \implies
      & \norma{y - x_{\bar h}} \leq r/2
    \end{array}
    \implies \norma{x-y} \leq r
    \implies \nabla \eta(x-y) = 0 \,.
  \end{displaymath}
  On the other hand, when $h \neq \bar h$ using~\ref{item:31}
  \begin{displaymath}
    \begin{array}{rcl}
      x \in B (x_{\bar h}, r/2)
      & \implies
      & \norma{x - x_{\bar h}} \leq r/2
      \\
      y \in B (x_h, r/2)
      & \implies
      & \norma{y - x_h} \leq r/2
    \end{array}
    \implies \norma{x-y} \geq \ell
    \implies \nabla \eta(x-y) = 0 \,.
  \end{displaymath}
  Thus, the velocity $v_i$ in~\eqref{eq:25} vanishes for all $t\geq 0$
  and $x \in \spt\rho_o$. Hence, $\rho_i (t,x) \, v_i (t,x)=0$ for all
  $t\geq 0$ and $x \in \reali^n$, completing the proof.
\end{proofof}

\section{Conclusion}

This paper presents an equation, namely~\eqref{eq:1}, establishes its
well posedness and proves some of its qualitative properties. The
resulting picture is that of a model prone to describe
\emph{collective behaviors} in the widest sense. The literature offers
many papers following the opposite direction, i.e., starting on the
basis of quantitative experiments, or on qualitative/heuristic
insights to specific real phenomena, and then elaborating up to getting
to a model. Many of the specific real features captured by these
models appear to be within the range of the qualitative properties of
the solutions to~\eqref{eq:1}.

The absence of a parabolic diffusion allows a (non standard) analytic
treatment essentially based on hyperbolic techniques. Moreover, a few
qualitative properties of the solutions to~\eqref{eq:1} are rigorously
proved, then exemplified through numerical integrations.

This new model suggests various related open analytical problems, such
as its connections with microscopic models, the characterization of
its asymptotic states, the development of \emph{ad hoc} efficient
numerical algorithms, for instance.

Moreover, the use of this equation as an encryption -- decryption tool
rises new analytical and numerical questions. The search for a class
where the map
$(\eta,V) \mapsto \mathcal{G}^{\eta,V}_{\strut T} \rho_o$ turns out to
be surjective is likely to require the development of tools for the
fine analysis of non local conservation laws. Present general
numerical methods for the integration of conservation laws hardly
respect the reversibility of~\eqref{eq:1}, in particular in more than
$1$ space dimension.

\appendix
\section{Time Reversibility}
\label{sec:time-reversibility}

In this appendix, for readability purposes, we insert the proof of
reversibility in time, which makes the map $\mathcal G$ in
\Cref{thm:cauchyProblem} a group.  More precisely, using the notation
in \textbf{Claim 5} in the proof of \Cref{thm:cauchyProblem}, we prove
that, for
$\rho_o \in \left(\L1 \cap \L\infty \cap \BV\right)(\reali^n;
\reali^m)$, the relation
$\mathcal G^-_t (\mathcal G^+_t \rho_o) = \rho_o$ holds for every
$t \in \reali_+$.

Fix
$\rho_o \in \left(\L1 \cap \L\infty \cap \BV\right)(\reali^n;
\reali^m)$, $t > 0$ and define, for $\tau \in [0, t]$,
$r(\tau) = \mathcal G^+_{t - \tau} \rho_o$. We prove that $r$
solves~\eqref{eq:48} with initial condition $\mathcal G^+_{t} \rho_o$,
so that $r(\tau) = \mathcal G^-_{\tau} (\mathcal G^+_t \rho_o)$.  To
this aim, fix $i \in \left\{1, \ldots, m\right\}$ and
$\psi \in \Cc\infty\left([0, t] \times \reali^n; \reali\right)$ and
use the change of variables $s = t - \tau$:
\begin{align*}
  &
    \int_0^t \int_{\reali^n} \left[r_i(\tau, x) \partial_t
    \psi(\tau, x) - r_i(\tau, x) V_i\left(\left(\nabla r(\tau) * \eta\right)(x)
    \right) \cdot \nabla \psi(\tau, x)\right] \d x \d \tau
  \\
  & \quad +\int_{\reali^n} r(0, x) \psi(0, x) \d x
    - \int_{\reali^n} r(t, x) \psi(t, x) \d x
  \\
  =
  & \int_0^t \int_{\reali^n} \left[(\mathcal G^+_{t - \tau} \rho_o)_i(x) \partial_t
    \psi(\tau, x) - (\mathcal G^+_{t - \tau} \rho_o)_i(x) V_i\left(\left(
    \nabla r(\tau) * \eta\right)(x)
    \right) \cdot \nabla \psi(\tau, x)\right] \d x \d \tau
  \\
  & \quad +\int_{\reali^n} \mathcal G^+_t \rho_o(x) \psi(0, x) \d x
    - \int_{\reali^n} \rho_o (x) \psi(t, x) \d x
  \\
  =
  & \int_0^t \int_{\reali^n} \left[(\mathcal G^+_{s} \rho_o)_i(x) \partial_t
    \psi(t - s, x) - (\mathcal G^+_{s} \rho_o)_i(x) V_i\left(\left(\nabla r
    (t-s) * \eta\right)(x)
    \right) \cdot \nabla \psi(t-s, x)\right] \d x \d s
  \\
  & \quad +\int_{\reali^n} \mathcal G^+_t \rho_o(x) \psi(0, x) \d x
    - \int_{\reali^n} \rho_o (x) \psi(t, x) \d x.
\end{align*}
Define, for every $s \in [0, t]$ and $x \in \reali^n$, the function
$\tilde \psi(s, x) = \psi(t-s, x)$, so that
$\tilde \psi \in \Cc\infty\left([0, t] \times \reali^n;
  \reali\right)$,
$\partial_s \tilde \psi (s, x) = - \partial_t \psi(t-s, x)$, and
$\nabla \tilde \psi (s, x) = \nabla \psi(t-s, x)$.  Thus
\begin{align*}
  &
    \int_0^t \int_{\reali^n} \left[r_i(\tau, x) \partial_t
    \psi(\tau, x) - r_i(\tau, x) V_i\left(\left(\nabla r(\tau) * \eta\right)(x)
    \right) \cdot \nabla \psi(\tau, x)\right] \d x \d \tau
  \\
  & \quad +\int_{\reali^n} r(0, x) \psi(0, x) \d x
    - \int_{\reali^n} r(t, x) \psi(t, x) \d x
  \\
  =
  & -\int_0^t \int_{\reali^n} \left[(\mathcal G^+_{s} \rho_o)_i(x) \partial_s
    \tilde \psi(s, x) + (\mathcal G^+_{s} \rho_o)_i(x)
    V_i\left(\left(\nabla r(t-s) * \eta\right)(x)
    \right) \cdot \nabla \tilde \psi(s, x)\right] \d x \d s
  \\
  & \quad +\int_{\reali^n} \mathcal G^+_t \rho_o(x) \tilde \psi(t, x) \d x
    - \int_{\reali^n} \rho_o (x) \tilde \psi(0, x) \d x
  \\
  =
  & -\int_0^t \int_{\reali^n} \left[(\mathcal G^+_{s} \rho_o)_i(x) \partial_s
    \tilde \psi(s, x) + (\mathcal G^+_{s} \rho_o)_i(x) V_i\left(\left(\nabla
    (\mathcal G^+_{s} \rho_o)* \eta\right)(x)
    \right) \cdot \nabla \tilde \psi(s, x)\right] \d x \d s
  \\
  & \quad +\int_{\reali^n} \mathcal G^+_t \rho_o(x) \tilde \psi(t, x) \d x
    - \int_{\reali^n} \rho_o (x) \tilde \psi(0, x) \d x
  = 0
\end{align*}
since $\tau \mapsto \mathcal G^+_\tau \rho_o$ solves the forward
problem. The uniqueness proved in Theorem~\ref{thm:cauchyProblem} thus
implies that $r(\tau) = \mathcal G^-_{\tau} (\mathcal G^+_t \rho_o)$,
and so $\mathcal G^-_t (\mathcal G^+_t \rho_o) = \rho_o$, completing
the proof.

\section*{Acknowledgments}

Both authors were partly supported by the GNAMPA~2022 project
\emph{Evolution Equations: Well Posedness, Control and Applications}
and by the PRIN~2022 project \emph{Modeling, Control and Games through
  Partial Differential Equations} (D53D23005620006), funded by the
European Union - Next Generation EU.

M.~Garavello was partially supported by the project funded under the
National Recovery and Resilience Plan (NRRP) of Italian Ministry of
University and Research funded by the European Union-NextGenerationEU.
Award Number: CN000023, Concession Decree No.  1033 of 17/06/2022
adopted by the Italian Ministry of University and Research, CUP:
H43C22000510001, Centro Nazionale per la Mobilità Sostenibile.

{\small

  \bibliographystyle{abbrv}

  \bibliography{GeneralClusters}

}

\end{document}